\documentclass[12pt,a4paper]{article}
\usepackage{amsmath,amsthm,amssymb,latexsym,epic,bezier}

\newtheorem{theorem}{Theorem}
\newtheorem{lemma}[theorem]{Lemma}
\newtheorem{corollary}[theorem]{Corollary}
\newtheorem{proposition}[theorem]{Proposition}

\newtheorem{remark}[theorem]{Remark}

\usepackage[all]{xy}
\usepackage{enumerate}
\numberwithin{equation}{section}

\newcommand{\IS}{\mathcal{IS}}

\renewcommand{\phi}{\varphi}
\newcommand{\tto}{\twoheadrightarrow}

\begin{document}

\title{On presentations of Brauer-type monoids}
\author{Ganna Kudryavtseva and Volodymyr Mazorchuk}
\date{}
\maketitle

\begin{abstract}
We obtain presentations for the Brauer monoid, the 
partial analogue of the Brauer monoid, and
for the greatest factorizable inverse submonoid
of the dual symmetric inverse monoid. In all three cases
we apply the same approach, based on the realization of all
these monoids as Brauer-type monoids. 
\end{abstract}

\section{Introduction and preliminaries}\label{s1}

The classical Coxeter presentation of the symmetric group $S_n$ plays
an important role in many branches of modern mathematics and physics. 
In the semigroup theory there are several ``natural'' analogues
of the symmetric group. For example the symmetric inverse
semigroup $\mathcal{IS}_n$ or the full transformation semigroup
$\mathcal{T}_n$. Perhaps a ``less natural'' generalization of
$S_n$ is the so-called Brauer semigroup $\mathfrak{B}_n$, which 
appeared in the context of centralizer algebras in representation 
theory in \cite{Br}. The basis of this algebra can be described in 
a nice combinatorial way using special {\em diagrams}
(see Section~\ref{s2}).  This combinatorial description motivated
a generalization of the Brauer algebra, the so-called {\em partition 
algebra}, which has its origins in physics and topology, see 
\cite{Mar}, \cite{Jo}. 
This algebra leads to another finite semigroup, the {\em partition 
semigroup}, usually denoted by $\mathfrak{C}_n$. Many classical 
semigroups, in particular, $S_n$, $\mathcal{IS}_n$, $\mathfrak{B}_n$ 
and some others (again see Section~\ref{s2}) are subsemigroups in
$\mathfrak{C}_n$. 

In the present paper we address the question of finding a presentation 
for some subsemigroups of $\mathfrak{C}_n$. As we have already 
mentioned, for $S_n$ this is a famous and very important result,
where the major role is played by the so-called {\em braid 
relations}. Because of the ``geometric'' nature of the generators
of the semigroups we consider, our initial motivation was that 
the additional relations for our semigroups would be some kind of
``singular deformations'' of the braid relations (analogous to 
the case of the singular braid monoid, see \cite{Ba, Bi}, or to
the known presentations of the Brauer algebra from 
\cite{BR}, \cite{BW}). In particular,  we wanted to get a complete 
list of ``deformations'' of the braid relations, which can appear 
in our cases. It turns out the all the semigroups we considered
indeed have presentations, all ingredients of which 
are in some sense deformations or degenerations of the braid 
relations.

As the main results of the paper we obtain a presentation for 
the semigroup $\mathfrak{B}_n$ (see Section~\ref{s3}), its
partial analogue $\mathcal{P}\mathfrak{B}_n$ (which can be also
called the {\em rook Brauer monoid}, see Section~\ref{s5}, and is a
kind of mixture of $\mathfrak{B}_n$ and $\mathcal{IS}_n$),
and a special inverse subsemigroup $\mathcal{IT}_n$ of 
$\mathfrak{C}_n$, which is isomorphic to the greatest 
factorizable inverse submonoid of the dual symmetric inverse 
monoid, see Section~\ref{s4} (another presentation for the
latter monoid was obtained in \cite{Fi}). 
The technical details in all cases are quite different,
however, the general approach is the same. We first ``guess''
the relations and in the standard way obtain an epimorphism
from the semigroup $T$, given by the corresponding presentation,
onto the semigroup we are dealing with. 
The only problem is to show that this epimorphism is in fact
a bijection. For this we have to compare the cardinalities of the
semigroups. In all our cases the symmetric group $S_n$ is the group 
of units in $T$. The product  $S_n\times S_n$ thus acts on $T$ 
via multiplication from the left and from the right. The idea is 
to show that  each orbit of this action contains a very special 
element, for which, using the relations, one can estimate the
cardinality of the stabilizer. The necessary statement then follows 
by comparing the cardinalities. 

\vspace{5mm}

\noindent
{\bf Acknowledgments.} The paper was written during the visit of the 
first author to Uppsala University, which was supported by the 
Swedish Institute. The financial support of the Swedish Institute and 
the hospitality of Uppsala University are gratefully acknowledged. For 
the second author the research was partially supported by the Swedish 
Research Council. We thank Victor Maltcev for informing us about the 
reference \cite{Fi}. We would also like to thank the referee for very
helpful suggestions.
\vspace{0.5cm}

\section{Brauer type semigroups}\label{s2}

For $n\in\mathbb{N}$ we denote by $S_n$ the {\em symmetric group} of
all permutations on the set $\{1,2,\dots,n\}$. We will consider 
the natural {\em right} action of $S_n$ on $\{1,2,\dots,n\}$ and the
induced action on the Boolean of $\{1,2,\dots,n\}$.
For a semigroup, $S$, we denote by $E(S)$ the set of all
idempotents of $S$.

Fix $n\in\mathbb{N}$ and let $M=M_n=\{1,2,\dots,n\}$,
$M'=\{1',2',\dots,n'\}$. We will consider ${}':M\to M'$ as a bijection,
whose inverse we will also denote by ${}'$.

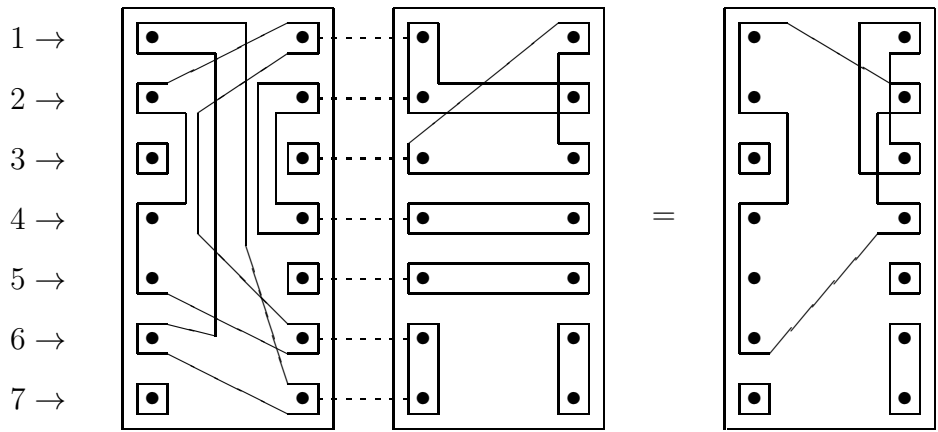
\begin{figure}
\special{em:linewidth 0.4pt}
\unitlength 0.80mm
\linethickness{0.4pt}
\begin{center}
\begin{picture}(150.00,85.00)
\drawline(15.00,10.00)(50.00,10.00)
\drawline(60.00,10.00)(95.00,10.00)
\drawline(115.50,10.00)(150.00,10.00)
\drawline(15.00,80.00)(50.00,80.00)
\drawline(60.00,80.00)(95.00,80.00)
\drawline(115.00,80.00)(150.00,80.00)
\drawline(15.00,10.00)(15.00,80.00)
\drawline(50.00,10.00)(50.00,80.00)
\drawline(60.00,10.00)(60.00,80.00)
\drawline(95.00,10.00)(95.00,80.00)
\drawline(115.00,10.00)(115.00,80.00)
\drawline(150.00,10.00)(150.00,80.00)
\put(1.00,75.00){\makebox(0,0)[cc]{$1\rightarrow$}}
\put(1.00,65.00){\makebox(0,0)[cc]{$2\rightarrow$}}
\put(1.00,55.00){\makebox(0,0)[cc]{$3\rightarrow$}}
\put(1.00,45.00){\makebox(0,0)[cc]{$4\rightarrow$}}
\put(1.00,35.00){\makebox(0,0)[cc]{$5\rightarrow$}}
\put(1.00,25.00){\makebox(0,0)[cc]{$6\rightarrow$}}
\put(1.00,15.00){\makebox(0,0)[cc]{$7\rightarrow$}}
\put(20.00,15.00){\makebox(0,0)[cc]{$\bullet$}}
\put(20.00,25.00){\makebox(0,0)[cc]{$\bullet$}}
\put(20.00,35.00){\makebox(0,0)[cc]{$\bullet$}}
\put(20.00,45.00){\makebox(0,0)[cc]{$\bullet$}}
\put(20.00,55.00){\makebox(0,0)[cc]{$\bullet$}}
\put(20.00,65.00){\makebox(0,0)[cc]{$\bullet$}}
\put(20.00,75.00){\makebox(0,0)[cc]{$\bullet$}}
\put(45.00,15.00){\makebox(0,0)[cc]{$\bullet$}}
\put(45.00,25.00){\makebox(0,0)[cc]{$\bullet$}}
\put(45.00,35.00){\makebox(0,0)[cc]{$\bullet$}}
\put(45.00,45.00){\makebox(0,0)[cc]{$\bullet$}}
\put(45.00,55.00){\makebox(0,0)[cc]{$\bullet$}}
\put(45.00,65.00){\makebox(0,0)[cc]{$\bullet$}}
\put(45.00,75.00){\makebox(0,0)[cc]{$\bullet$}}
\put(65.00,15.00){\makebox(0,0)[cc]{$\bullet$}}
\put(65.00,25.00){\makebox(0,0)[cc]{$\bullet$}}
\put(65.00,35.00){\makebox(0,0)[cc]{$\bullet$}}
\put(65.00,45.00){\makebox(0,0)[cc]{$\bullet$}}
\put(65.00,55.00){\makebox(0,0)[cc]{$\bullet$}}
\put(65.00,65.00){\makebox(0,0)[cc]{$\bullet$}}
\put(65.00,75.00){\makebox(0,0)[cc]{$\bullet$}}
\put(90.00,15.00){\makebox(0,0)[cc]{$\bullet$}}
\put(90.00,25.00){\makebox(0,0)[cc]{$\bullet$}}
\put(90.00,35.00){\makebox(0,0)[cc]{$\bullet$}}
\put(90.00,45.00){\makebox(0,0)[cc]{$\bullet$}}
\put(90.00,55.00){\makebox(0,0)[cc]{$\bullet$}}
\put(90.00,65.00){\makebox(0,0)[cc]{$\bullet$}}
\put(90.00,75.00){\makebox(0,0)[cc]{$\bullet$}}
\put(120.00,15.00){\makebox(0,0)[cc]{$\bullet$}}
\put(120.00,25.00){\makebox(0,0)[cc]{$\bullet$}}
\put(120.00,35.00){\makebox(0,0)[cc]{$\bullet$}}
\put(120.00,45.00){\makebox(0,0)[cc]{$\bullet$}}
\put(120.00,55.00){\makebox(0,0)[cc]{$\bullet$}}
\put(120.00,65.00){\makebox(0,0)[cc]{$\bullet$}}
\put(120.00,75.00){\makebox(0,0)[cc]{$\bullet$}}
\put(145.00,15.00){\makebox(0,0)[cc]{$\bullet$}}
\put(145.00,25.00){\makebox(0,0)[cc]{$\bullet$}}
\put(145.00,35.00){\makebox(0,0)[cc]{$\bullet$}}
\put(145.00,45.00){\makebox(0,0)[cc]{$\bullet$}}
\put(145.00,55.00){\makebox(0,0)[cc]{$\bullet$}}
\put(145.00,65.00){\makebox(0,0)[cc]{$\bullet$}}
\put(145.00,75.00){\makebox(0,0)[cc]{$\bullet$}}
\put(105.00,45.00){\makebox(0,0)[cc]{$=$}}
\dashline{1}(45.00,75.00)(65.00,75.00)
\dashline{1}(45.00,65.00)(65.00,65.00)
\dashline{1}(45.00,55.00)(65.00,55.00)
\dashline{1}(45.00,45.00)(65.00,45.00)
\dashline{1}(45.00,35.00)(65.00,35.00)
\dashline{1}(45.00,25.00)(65.00,25.00)
\dashline{1}(45.00,15.00)(65.00,15.00)
\linethickness{0.2pt}
\drawline(17.50,12.50)(22.50,12.50)
\drawline(17.50,17.50)(22.50,17.50)
\drawline(17.50,12.50)(17.50,17.50)
\drawline(22.50,12.50)(22.50,17.50)
\drawline(17.50,52.50)(22.50,52.50)
\drawline(17.50,57.50)(22.50,57.50)
\drawline(17.50,52.50)(17.50,57.50)
\drawline(22.50,52.50)(22.50,57.50)
\drawline(117.50,12.50)(122.50,12.50)
\drawline(117.50,17.50)(122.50,17.50)
\drawline(117.50,12.50)(117.50,17.50)
\drawline(122.50,12.50)(122.50,17.50)
\drawline(117.50,52.50)(122.50,52.50)
\drawline(117.50,57.50)(122.50,57.50)
\drawline(117.50,52.50)(117.50,57.50)
\drawline(122.50,52.50)(122.50,57.50)
\drawline(47.50,32.50)(42.50,32.50)
\drawline(47.50,37.50)(42.50,37.50)
\drawline(47.50,32.50)(47.50,37.50)
\drawline(42.50,32.50)(42.50,37.50)
\drawline(47.50,52.50)(42.50,52.50)
\drawline(47.50,57.50)(42.50,57.50)
\drawline(47.50,52.50)(47.50,57.50)
\drawline(42.50,52.50)(42.50,57.50)
\drawline(147.50,32.50)(142.50,32.50)
\drawline(147.50,37.50)(142.50,37.50)
\drawline(147.50,32.50)(147.50,37.50)
\drawline(142.50,32.50)(142.50,37.50)
\drawline(67.50,12.50)(62.50,12.50)
\drawline(67.50,27.50)(62.50,27.50)
\drawline(67.50,12.50)(67.50,27.50)
\drawline(62.50,12.50)(62.50,27.50)
\drawline(87.50,12.50)(92.50,12.50)
\drawline(87.50,27.50)(92.50,27.50)
\drawline(87.50,12.50)(87.50,27.50)
\drawline(92.50,12.50)(92.50,27.50)
\drawline(147.50,12.50)(142.50,12.50)
\drawline(147.50,27.50)(142.50,27.50)
\drawline(147.50,12.50)(147.50,27.50)
\drawline(142.50,12.50)(142.50,27.50)
\drawline(117.50,22.50)(117.50,47.50)
\drawline(125.50,47.50)(117.50,47.50)
\drawline(125.50,47.50)(125.50,62.50)
\drawline(117.50,62.50)(125.50,62.50)
\drawline(117.50,62.50)(117.50,77.50)
\drawline(125.50,77.50)(117.50,77.50)
\drawline(125.50,77.50)(142.50,67.50)
\drawline(147.50,67.50)(142.50,67.50)
\drawline(147.50,67.50)(147.50,62.50)
\drawline(140.50,62.50)(147.50,62.50)
\drawline(140.50,62.50)(140.50,47.50)
\drawline(147.50,47.50)(140.50,47.50)
\drawline(147.50,47.50)(147.50,42.50)
\drawline(140.50,42.50)(147.50,42.50)
\drawline(140.50,42.50)(122.50,22.50)
\drawline(122.50,22.50)(117.50,22.50)
\drawline(147.50,52.50)(147.50,57.50)
\drawline(142.50,57.50)(147.50,57.50)
\drawline(142.50,57.50)(142.50,72.50)
\drawline(147.50,72.50)(142.50,72.50)
\drawline(147.50,72.50)(147.50,77.50)
\drawline(137.50,77.50)(147.50,77.50)
\drawline(137.50,77.50)(137.50,52.50)
\drawline(147.50,52.50)(137.50,52.50)
\drawline(62.50,32.50)(92.50,32.50)
\drawline(62.50,37.50)(92.50,37.50)
\drawline(62.50,42.50)(92.50,42.50)
\drawline(62.50,47.50)(92.50,47.50)
\drawline(62.50,52.50)(92.50,52.50)
\drawline(62.50,62.50)(92.50,62.50)
\drawline(62.50,32.50)(62.50,37.50)
\drawline(92.50,32.50)(92.50,37.50)
\drawline(62.50,42.50)(62.50,47.50)
\drawline(92.50,42.50)(92.50,47.50)
\drawline(62.50,52.50)(62.50,57.50)
\drawline(92.50,52.50)(92.50,57.50)
\drawline(92.50,62.50)(92.50,67.50)
\drawline(92.50,72.50)(92.50,77.50)
\drawline(67.50,67.50)(92.50,67.50)
\drawline(67.50,67.50)(67.50,77.50)
\drawline(62.50,77.50)(67.50,77.50)
\drawline(62.50,77.50)(62.50,62.50)
\drawline(62.50,57.50)(87.50,77.50)
\drawline(92.50,77.50)(87.50,77.50)
\drawline(92.50,72.50)(87.50,72.50)
\drawline(87.50,57.50)(87.50,72.50)
\drawline(87.50,57.50)(92.50,57.50)
\drawline(47.50,42.50)(47.50,47.50)
\drawline(40.50,47.50)(47.50,47.50)
\drawline(40.50,47.50)(40.50,62.50)
\drawline(47.50,62.50)(40.50,62.50)
\drawline(47.50,62.50)(47.50,67.50)
\drawline(37.50,67.50)(47.50,67.50)
\drawline(37.50,67.50)(37.50,42.50)
\drawline(47.50,42.50)(37.50,42.50)
\drawline(17.50,32.50)(17.50,47.50)
\drawline(25.50,47.50)(17.50,47.50)
\drawline(25.50,47.50)(25.50,62.50)
\drawline(17.50,62.50)(25.50,62.50)
\drawline(17.50,62.50)(17.50,67.50)
\drawline(22.50,67.50)(17.50,67.50)
\drawline(22.50,67.50)(42.50,77.50)
\drawline(47.50,77.50)(42.50,77.50)
\drawline(47.50,77.50)(47.50,72.50)
\drawline(42.50,72.50)(47.50,72.50)
\drawline(17.50,32.50)(22.50,32.50)
\drawline(42.50,22.50)(22.50,32.50)
\drawline(42.50,22.50)(47.50,22.50)
\drawline(47.50,27.50)(47.50,22.50)
\drawline(47.50,27.50)(42.50,27.50)
\drawline(27.50,42.50)(42.50,27.50)
\drawline(27.50,42.50)(27.50,62.50)
\drawline(42.50,72.50)(27.50,62.50)
\drawline(17.50,22.50)(17.50,27.50)
\drawline(17.50,22.50)(22.50,22.50)
\drawline(42.50,12.50)(22.50,22.50)
\drawline(42.50,12.50)(47.50,12.50)
\drawline(47.50,17.50)(47.50,12.50)
\drawline(47.50,17.50)(42.50,17.50)
\drawline(17.50,72.50)(17.50,77.50)
\drawline(35.50,77.50)(17.50,77.50)
\drawline(35.50,77.50)(35.50,40.50)
\drawline(42.50,17.50)(35.50,40.50)
\drawline(17.50,72.50)(30.50,72.50)
\drawline(30.50,25.50)(30.50,72.50)
\drawline(30.50,25.50)(22.50,27.50)
\drawline(17.50,27.50)(22.50,27.50)
\end{picture}
\end{center}
\caption{Multiplication of elements of $\mathfrak{C}_n$.}\label{fig1}
\end{figure}

Consider the set $\mathfrak{C}_n$ 
of all decompositions of $M\cup M'$ into disjoint unions of subsets. 
Given $\alpha,\beta\in \mathfrak{C}_n$, $\alpha=X_1\cup\dots \cup X_k$
and $\beta=Y_1\cup\dots \cup Y_l$, we define their product
$\gamma=\alpha\beta$ as the unique element of $\mathfrak{C}_n$
satisfying the following conditions:
\begin{enumerate}[(P1)]
\item\label{mult1} For $i,j\in M$ the elements $i$ and $j$ belong
to the same block  of the decomposition $\gamma$ if an only if
they belong to the same block of the decomposition $\alpha$ or there
exists a sequence, $s_1,\dots,s_{m}$, where $m$ is even, 
of elements from $M$ such that
$i$ and $s'_1$  belong to the same block of $\alpha$;
$s_1$ and $s_2$  belong to the same block of $\beta$;
$s'_2$ and $s'_3$  belong to the same block of $\alpha$ and so on;
$s_{m-1}$ and $s_{m}$  belong to the same block of $\beta$;
$s'_{m}$ and $j$  belong to the same block of $\alpha$.
\item\label{mult2} For $i,j\in M$ the elements $i'$ and $j'$ belong
to the same block  of the decomposition $\gamma$ if an only if
they belong to the same block of the decomposition $\beta$ or there
exists a sequence, $s_1,\dots,s_{m}$, where $m$ is even, 
of elements from $M$ such that
$i'$ and $s_1$  belong to the same block of $\beta$;
$s'_1$ and $s'_2$  belong to the same block of $\alpha$;
$s_2$ and $s_3$  belong to the same block of $\beta$ ans so on;
$s'_{m-1}$ and $s'_{m}$  belong to the same block of $\alpha$;
$s_{m}$ and $j'$  belong to the same block of $\beta$.
\item\label{mult3} For $i,j\in M$ the elements $i$ and $j'$ belong
to the same block  of the decomposition $\gamma$ if an only if
there exists a sequence, $s_1,\dots,s_{m}$, where $m$ is odd,
of elements from $M$ such that $i$ and $s'_1$  belong to the 
same block of $\alpha$; $s_1$ and $s_2$  belong to the same block of 
$\beta$; $s'_2$ and $s'_3$  belong to the same block of $\alpha$ and so on;
$s'_{m-1}$ and $s'_{m}$  belong to the same block of $\alpha$;
$s_{m}$ and $j'$  belong to the same block of $\beta$.
\end{enumerate}
One can think about the elements of $\mathfrak{C}_n$ as ``microchips''
or ``generalized microchips'' with $n$ pins on the left hand side
(corresponding to the elements of $M$) and $n$ pins on the right 
hand side (corresponding to the elements of $M'$). For 
$\alpha\in \mathfrak{C}_n$ we connect two pins of the corresponding
chip if and only if they belong to the same set of the partition
$\alpha$. The operation described above can then be viewed as
a ``composition'' of such chips: having $\alpha,\beta\in\mathfrak{C}_n$ 
we identify (connect) the right pins of $\alpha$ with the corresponding 
left pins of $\beta$, which uniquely defines a connection of the remaining 
pins (which are the left pins of $\alpha$ and the right pins of $\beta$). 
An example of multiplication of two chips from $\mathfrak{C}_n$ is given on Figure~\ref{fig1}. Note that, performing the operation we can obtain 
some ``dead circles'' formed by some identified pins from $\alpha$ and 
$\beta$. These circles should be disregarded (however they play an
important role in representation theory as they allow to deform the
multiplication in the semigroup algebra). From this interpretation
it is fairly obvious  that the composition of elements from $\mathfrak{C}_n$ 
defined above is associative. On the level of associative algebra,
the partition algebra was defined in \cite{Mar} and then studied by
several authors especially in recent years, see for example 
\cite{Bl,Mar2,MarEl,MarWo,Pa,Xi}. Purely
as a semigroup it seems that $\mathfrak{C}_n$ appeared in \cite{Maz}.

Let $\alpha\in \mathfrak{C}_n$ and $X$ be a block of $\alpha$. The block 
$X$ will be called 
\begin{itemize}
\item a {\em line} provided that $\vert X\vert=2$ and $X$ intersects 
with both $M$ and  $M'$;
\item a {\em generalized line} provided that  $X$ intersects with 
both $M$ and $M'$;
\item a {\em bracket} if $\vert X\vert=2$ and either $X\subset M$ or 
$X\subset M'$;
\item a {\em generalized bracket} if $\vert X\vert\geq 2$ and either 
$X\subset M$ or $X\subset M'$;
\item a {\em point} if $\vert X\vert=1$.
\end{itemize}
By a {\em Brauer-type semigroup} we will mean a ``natural'' subsemigroup 
of the semigroup $\mathfrak{C}_n$. Here are some examples:
\begin{enumerate}[(E1)]
\item  The subsemigroup, consisting of all elements 
$\alpha\in  \mathfrak{C}_n$ such that each block of $\alpha$ is a line. 
This subsemigroup is canonically identified with $S_n$ and is the 
group of units of $\mathfrak{C}_n$.
\item The subsemigroup, consisting of all elements 
$\alpha\in  \mathfrak{C}_n$ such that each block of
$\alpha$ is a either a line or a point. This subsemigroup is 
canonically identified with the {\em symmetric inverse
semigroup} $\mathcal{IS}_n$.
\item The subsemigroup $\mathfrak{B}_n$, consisting of all elements 
$\alpha\in  \mathfrak{C}_n$ such that each block of
$\alpha$ is a either a line or a bracket. This is the classical
{\em Brauer semigroup}, see \cite{Ke,Maz0}.
\item The subsemigroup $\mathcal{P}\mathfrak{B}_n$, consisting of 
all elements $\alpha\in  \mathfrak{C}_n$ such that each block of
$\alpha$ is a either a line or a bracket or a point. This is the 
{\em partial analogue of the Brauer semigroup}, see \cite{Maz0}.
\item The subsemigroup $\mathcal{IP}_n$, consisting of all 
$\alpha\in \mathfrak{C}_n$ such that each block of $\alpha$ is 
a generalized line. In this form the semigroup $\mathcal{IP}_n$
appeared in \cite{Mal1,Mal2}. It is easy to see that the 
semigroup $\mathcal{IP}_n$ is isomorphic to the {\em dual symmetric 
inverse monoid} $\mathcal{I}^*_M$ from \cite{FL}.
\item The subsemigroup $\mathcal{IT}_n$, consisting of all 
$\alpha\in \mathfrak{C}_n$ such that each block $X$ of $\alpha$ is  
a generalized line and $\vert X\cap M\vert=\vert X\cap M'\vert$.
In this form the semigroup $\mathcal{IT}_n$ appeared in \cite{Mal2}.
The semigroup $\mathcal{IT}_n$ is isomorphic to the 
{\em greatest factorizable inverse submonoid} $\mathcal{F}^*_M$ 
of $\mathcal{I}^*_M$ from \cite{FL}.
\end{enumerate}

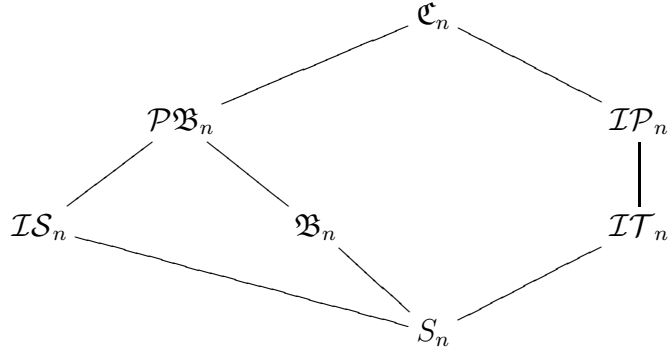
\begin{figure}
\special{em:linewidth 0.4pt}
\unitlength 0.80mm
\linethickness{0.4pt}
\begin{center}
\begin{displaymath}
\xymatrix{
 &&&\mathfrak{C}_n\ar@{-}[rrd]\ar@{-}[lld]&&  \\
 &\mathcal{P}\mathfrak{B}_n\ar@{-}[ld]\ar@{-}[rd]
 &&&&\mathcal{IP}_n\ar@{-}[d] \\
\mathcal{IS}_n&&\mathfrak{B}_n&&&\mathcal{IT}_n\\
&&&S_n\ar@{-}[rru]\ar@{-}[lu]\ar@{-}[lllu]&&
}
\end{displaymath}
\end{center}
\caption{Inclusions for classical Brauer-type semigroups}\label{fig2}
\end{figure}

All the semigroups described above are regular. $S_n$ is a group.
The semigroups $IS_n$, $\mathcal{IP}_n$ and $\mathcal{IT}_n$ are inverse, 
while $\mathfrak{C}_n$, $\mathfrak{B}_n$ and $\mathcal{P}\mathfrak{B}_n$ 
are not. The partially ordered set consisting of these semigroups, with
the partial order given by inclusions, is illustrated on
Figure~\ref{fig2}. 

In what follows we will need some easy combinatorial results
for Brauer-type semigroups. For $\alpha\in\mathfrak{C}_n$ we
define the {\em rank} $\mathrm{rk}(\alpha)$ of $\alpha$ as the number
of generalized lines in $\alpha$, that is the number of blocks in 
$\alpha$ intersecting with both $M$ and $M'$. Note that for the semigroups
$S_n$, $\mathcal{IS}_n$, $\mathfrak{B}_n$, $\mathcal{P}\mathfrak{B}_n$
and $\mathfrak{C}_n$ ranks of the elements classify the 
$\mathcal{D}$-classes (this is obvious for $S_n$, for 
$\mathcal{IS}_n$ this is an easy exercise, for $\mathfrak{B}_n$ 
and $\mathcal{P}\mathfrak{B}_n$ this can be found in \cite{Maz0}, 
and for  $\mathfrak{C}_n$ it can be obtained by arguments similar to
those from \cite{Maz0} for $\mathfrak{B}_n$).

For the semigroup $\mathcal{IT}_n$ we will need a different notion. 
Let $X$ be a finite set and $X=\cup_{i=1}^k X_k$  be a decomposition of $X$ 
into a union of pairwise disjoint subsets. For each $i$, $1\leq i\leq |X|$, 
let $m_i$ denote the number of subsets of this decomposition, whose 
cardinality equals $i$. The tuple $(m_1,\dots, m_{|X|})$ will be called 
the {\em type}  of the decomposition. 
Consider an element, $\alpha\in {\mathcal{IT}_n}$. 
By definition $\alpha$ is a decomposition of $M\cup M'$ into a disjoint 
union of subsets, whose intersections with $M$ and $M'$ have the
same cardinality. Let $(m_1,\dots, m_{2n})$ be the type of this 
decompositions (note that $m_i\neq 0$ only if $i$ is even). 
The element $\alpha$ induces  a decomposition of $M$ into disjoint 
subsets, whose blocks are intersections 
of the blocks of $\alpha$ with $M$. By the {\em type} of $\alpha$ we will 
mean the type of this decomposition of $M$, which is obviously equal to $(m_2,m_4,\dots,m_{2n})$.  The types of elements from
$\mathcal{IT}_n$ correspond bijectively to partitions of $n$
(a {\em partition}, $\lambda\vdash n$, of $n$ is a tuple, $\lambda=(\lambda_1,\dots,\lambda_k)$,
of positive integers such that $\lambda_1\geq\lambda_2\geq
\dots\geq\lambda_k$ and $\lambda_1+\cdots+\lambda_k=n$).  The types of 
the elements classify the  $\mathcal{D}$-classes  in $\mathcal{IT}_n$, see
\cite[Section~3]{FL}.

For the semigroup $\mathcal{P}\mathfrak{B}_n$ we will need a more
complicated technical tool. Although $\mathcal{D}$-classes are 
classified by ranks we will need to distinguish elements of a given 
rank, so we introduce the notion of a type.
For $\alpha\in \mathcal{P}\mathfrak{B}_n$
let $r$ denote the number of lines in $\alpha$; 
$b_1$ the number of brackets in $\alpha$, contained in $M$;
$b_2$ the number of brackets in $\alpha$, contained in $M'$;
$p_1$ the number of points in $\alpha$, contained in $M$;
$p_2$ the number of points in $\alpha$, contained in $M'$.
Obviously $n=r+2b_1+p_1=r+2b_2+p_2$. Define the {\em type} of
$\alpha$ as follows:
\begin{displaymath}
\mathrm{type}(\alpha)=
\begin{cases}
(b_2,b_1-b_2,0,p_1), & b_1\geq b_2;\\
(b_1,0,b_2-b_1,p_2), & b_2>b_1.
\end{cases}
\end{displaymath}

We will need the following explicit combinatorial formulae for 
the number of elements of a  given rank or type.

\begin{proposition}\label{comb}
\begin{enumerate}[(a)]
\item\label{comb.1} For $k\in\{0,\dots,n\}$ the number of elements of
rank $k$ in $\mathcal{IS}_n$ equals $\binom{n}{k}^2 k!$.
\item\label{comb.2} For $k\in\{1,\dots,n\}$ the number of elements of
rank $k$ in $\mathfrak{B}_n$ equals $0$ if $n-k$ is odd and
$\frac{(n!)^2}{2^{2l}(l!)^2k!}$ if $n-k=2l$ is even.
\item\label{comb.3} The number of elements of $\mathcal{IT}_n$ of 
type $(m_1,\dots,m_n)$ equals 
\begin{displaymath}
\frac{(n!)^2} 
{\displaystyle\prod_{i=1}^n(m_i!(i!)^{2m_i})}.
\end{displaymath}
\item\label{comb.4} For all non-negative integers $k,m,t$ such that 
$2k+2m+t\leq n$ the number of elements of the type
$(k,m,0,t)$ in $\mathcal{P}\mathfrak{B}_n$ is equal to the
number of elements of the type
$(k,0,m,t)$ in $\mathcal{P}\mathfrak{B}_n$ and equals
\begin{displaymath}
\frac{(n!)^2}{k!2^k(t+2m)!(k+m)!2^{k+m}t!(n-2k-2m-t)!}.
\end{displaymath}
\end{enumerate}
\end{proposition}

\begin{proof}
This is a straightforward combinatorial calculation.
\end{proof}

\begin{remark}\label{remver}
{\rm
The semigroup $\mathfrak{C}_n$ can be also connected to some other
semigroups of binary relations. As we have already mentioned, the 
subsemigroup $\mathcal{IP}_n$ of $\mathfrak{C}_n$ is isomorphic to the
dual symmetric inverse monoid $\mathcal{I}^*_M$ from \cite{FL}, 
which is the semigroup of all difunctional binary relations under 
the operation of taking the smallest difunctional binary relations, 
containing the product of two given relations. The semigroup
$\mathcal{IT}_n$ is isomorphic to the greatest factorizable
inverse submonoid of $\mathcal{I}^*_M$, that is to the semigroup
$E(\mathcal{I}^*_M)S_n$. One can also deform the multiplication
in $\mathfrak{C}_n$ in the following way: given $\alpha,\beta\in
\mathfrak{C}_n$ define $\gamma=\alpha\star\beta$ as follows:
all blocks of $\gamma$ are either points or generalized lines, and
for $i,j\in M$ the elements $i$ and $j'$ belong to the same block
of $\gamma$ if and only if $i$ belongs to some block $X$ of $\alpha$
and $j'$ belongs to some block $Y$ of $\beta$ such that 
$X\cap M'=(Y\cap M)'$. It is straightforward that this deformed
multiplication is associative and hence we get a new semigroup,
$\tilde{\mathfrak{C}}_n$. This semigroup is an inflation of
Vernitsky's inverse semigroup $(D_X,\diamond)$, see \cite{Ve},
which is a subsemigroup of $\tilde{\mathfrak{C}}_n$ in the
natural way. An isomorphic object can be obtained if instead of
points one requires that $\gamma$ contains at most one generalized 
bracket,  which is a subset of $M$, and at most one generalized 
bracket, which  is a subset of $M'$.
}
\end{remark}

\section{Presentation for $\mathfrak{B}_n$}\label{s3}

For $i=1,\dots,n-1$ we denote by $s_i$ the elementary transposition
$(i,i+1)\in S_n$, and by $\pi_i$ the element
$\{i,i+1\}\cup \{i',(i+1)'\}\cup\bigcup_{j\neq i,i+1}\{j,j'\}$ of
$\mathfrak{B}_n$ (the elementary {\em atom} from \cite{Maz0}). 
It is easy to see (and can be derived from the results of
\cite{Maz0} and \cite{Mal0}) that $\mathfrak{B}_n$ is generated by 
$\{s_i\}\cup \{\pi_i\}$ as a monoid. Moreover, $\mathfrak{B}_n$ is even 
generated by $\{s_i\}$ and, for example, $\pi_1$. However, 
we think that the set $\{s_i\}\cup \{\pi_i\}$ is more natural as 
a system of generators for $\mathfrak{B}_n$, for example because of 
the connection between Brauer and Temperley-Lieb algebras (and analogy
with the singular braid monoid, see \cite{Ba,Bi}). In this  section 
we obtain a  presentation for $\mathfrak{B}_n$ with respect to this 
system of generators (this resembles the presentation of the
Brauer algebra in \cite{BW}, see also \cite{BR}).

Let $T$ denote the monoid with the identity element $e$, 
generated by the elements $\sigma_i$, $\theta_i$, 
$i=1,\dots,n-1$, subject to the following relations
(where $i,j\in\{1,2,\dots,n-1\}$):
\begin{gather}
\sigma_i^2=e; \quad
\sigma_i\sigma_j=\sigma_j\sigma_i,\,\, |i-j|>1;\quad
\sigma_i\sigma_j\sigma_i=\sigma_j\sigma_i\sigma_j, \,\, |i-j|=1;\label{it11}\\
\theta_i^2=\theta_i;\quad
\theta_i\theta_j=\theta_j\theta_i, \,\, |i-j|>1;\quad\theta_i\theta_j\theta_i=\theta_i,\,\, |i-j|=1;\label{it21}\\
\theta_i\sigma_i=\sigma_i\theta_i=\theta_i, \quad
\theta_i\sigma_j=\sigma_j\theta_i, \,\, |i-j|>1;\label{it31}\\
\sigma_{i}\theta_j\theta_i=\sigma_{j}\theta_i,\quad
\theta_i\theta_j\sigma_{i}=\theta_i\sigma_j,
\,\, |i-j|=1. \label{it41new}
\end{gather}

\begin{theorem}\label{theorem1}
The map $\sigma_i\mapsto s_i$ and $\theta_i\to \pi_i$, $i=1,\dots,n-1$,
extends to an isomorphism, $\varphi:T\to \mathfrak{B}_n$.
\end{theorem}

The rest of the section will be devoted to the proof of
Theorem~\ref{theorem1}. We start with the following easy observation,
which later on will be used in our computations:

\begin{lemma}\label{newrellem}
Under the assumption that the relations \eqref{it11}--\eqref{it41new}
are satisfied, we have the following relations:
\begin{gather}
\sigma_i\theta_j\sigma_i=\sigma_j\theta_i\sigma_j,\quad 
\theta_i\sigma_j\theta_i=\theta_i,\,\, |i-j|=1; \label{it41}\\
\sigma_i\sigma_{i+1}\theta_i\theta_{i+2}=
\sigma_{i+2}\sigma_{i+1}\theta_i\theta_{i+2}.\label{it51}
\end{gather}
\end{lemma}

\begin{proof}
For $i,j$, $|i-j|=1$, applying  \eqref{it41new} twice  we have
\begin{displaymath}
\sigma_i\theta_j\sigma_i=
\sigma_j\theta_i\theta_j\sigma_i=
\sigma_j\theta_i\sigma_j.
\end{displaymath}
Applying \eqref{it41new}, \eqref{it31} and, finally, \eqref{it21} 
we also have
\begin{displaymath}
\theta_i\sigma_j\theta_i=
\theta_i\theta_j\sigma_i\theta_i=
\theta_i\theta_j\theta_i=\theta_i.
\end{displaymath}
This gives \eqref{it41}. Analogously, applying
\eqref{it41new}, \eqref{it11}, \eqref{it21}  and
\eqref{it41new} again gives
\begin{displaymath}
\sigma_{i+2}\sigma_{i+1}\theta_i\theta_{i+2}=
\sigma_{i+2}\sigma_{i}\theta_{i+1}\theta_i\theta_{i+2}=
\sigma_{i}\sigma_{i+2}\theta_{i+1}\theta_{i+2}\theta_{i}=
\sigma_{i}\sigma_{i+1}\theta_{i+2}\theta_{i},
\end{displaymath}
which implies \eqref{it51}.
\end{proof}

It is a direct calculation to verify that the generators
$s_i$ and $\pi_i$ of $\mathfrak{B}_n$ satisfy the relations,
corresponding to \eqref{it11}--\eqref{it41new}. Thus the map
$\sigma_i\mapsto s_i$ and $\theta_i\mapsto \pi_i$, $i=1,\dots,n-1$,
extends to an epimorphism, $\varphi:T\tto \mathfrak{B}_n$. Hence,
to prove Theorem~\ref{theorem1} we have only to show that 
$|T|=|\mathfrak{B}_n|$. To do this we will have to study the
structure of the semigroup $T$ in details.

Let $W$ denote the free monoid, generated by $\sigma_i$,
$\theta_i$, $i=1,\dots,n-1$, and $\psi:W\tto T$ denote the
canonical projection. Let $\sim$ be the corresponding congruence on
$W$, that is $v\sim w$ provided that $\psi(v)=\psi(w)$. We
start with the following description of units in $T$:

\begin{lemma}\label{lemma301}
The elements $\sigma_i$, $i=1,\dots,n-1$, generate the group 
$G$ of units in $T$, which is isomorphic to the symmetric group $S_n$.
\end{lemma}

\begin{proof}
Let $v,w\in W$ be such that $v\sim w$. Assume further that
$v$ contains some $\theta_i$. Since $\theta$'s allways occur 
on both sides in the relations \eqref{it21}--\eqref{it41new} and 
do not occur in the relations \eqref{it11}, it follows that
$w$ must contain some $\theta_j$. In particular, the submonoid,
generated in $W$ by $\sigma_i$, $i=1,\dots,n-1$, is a union
of equivalence classes with respect to $\sim$. Using the well-known
Coxeter presentation of the symmetric group we obtain that 
$\sigma_i$, $i=1,\dots,n-1$, generate in $T$ a copy of the
symmetric group. All elements of this group are obviously units in $T$.
On the other hand, if $v,w\in W$ and $v$ contains some $\theta_i$,
then $vw$ contains $\theta_i$ as well. By the above arguments, 
$vw$ can not be equivalent to the empty word. Hence
$v$ is not invertable in $T$. The claim of the lemma follows.
\end{proof}

In what follows we will identify the group $G$ of units in $T$ with
$S_n$ via the isomorphism, which sends $\sigma_i\in G$ to $s_i$.
There is a natural action of $S_n$ on $T$ by inner automorphisms of 
$T$ via conjugation: $x^g=g^{-1}xg$ for each $x\in T$, $g\in S_n$. 

\begin{lemma}\label{lemma201}
The $S_n$-stabilizer of $\theta_1$ is the subgroup $H$ of
$S_n$, consisting of all permutations, which preserve the set
$\{1,2\}$. This subgroup is isomorphic to $S_2\times S_{n-2}$.
\end{lemma}

\begin{proof}
We have $\sigma_j\theta_1\sigma_j=\theta_j$, $j\neq 2$, by \eqref{it31}.
Since $\sigma_j$,  $j\neq 2$, generate $H$, we obtain that all
elements of $H$ stabilize $\theta_1$. In particular, the
$S_n$-orbit of $\theta_1$ consists of at most $|S_n|/|H|=\binom{n}{2}$
elements. At the same time, it is easy to see that the
$S_n$-orbit of $\varphi(\theta_1)$ consists of exactly
$\binom{n}{2}$ different elements and hence $H$ must coincide
with the $S_n$-stabilizer of $\theta_1$.
\end{proof}

Since $S_n$ acts on $T$ via automorphisms and $\theta_1$ is an idempotent,
all elements in the $S_n$-orbit of $\theta_1$ are idempotents.
From Lemma~\ref{lemma201} it follows that the elements of the
$S_n$-orbit of $\theta_1$ are in the natural bijection with the
cosets $H\backslash S_n$. By the definition of $H$, two elements,
$x,y\in S_n$, are contained in the same coset if and only if 
$x(\{1,2\})=y(\{1,2\})$.

\begin{lemma}\label{lemma202}
The $S_n$-orbit of $\theta_1$ contains all $\theta_i$,
$i=1,\dots,n-1$. Moreover, for $w\in S_n$ we have
$w^{-1}\theta_1 w=\theta_i$ if and only if $w(\{1,2\})=\{i,i+1\}$.
\end{lemma}

\begin{proof}
We use induction on $i$ with the case $i=1$ being trivial.
Let $i>1$ and assume that $\theta_{i-1}$ is contained in our orbit.
Then $\theta_{i}=\sigma_{i-1}\sigma_{i}\theta_{i-1}\sigma_{i}\sigma_{i-1}$
and hence $\theta_{i}$ is contained in our orbit as well. Hence
all $\theta_i$ indeed belong to the $S_n$-orbit of $\theta_1$.
The second claim follows from 
\begin{equation}\label{oldequation}
\sigma_{i-1}\sigma_{i}\sigma_{i-2}\sigma_{i-1}\cdots
\sigma_{1}\sigma_{2}(\{1,2\})=\{i,i+1\},
\end{equation}
which is obtained by a direct calculation. This completes the proof.
\end{proof}

For $w\in S_n$ such that $w(\{1,2\})=\{i,j\}$, where $i<j$, we
set $\epsilon_{i,j}=w^{-1}\theta_1 w$, which is well defined
by Lemma~\ref{lemma201}.

\begin{lemma}\label{lemma203}
Suppose $\{i,j\}\cap\{p,q\}=\varnothing$. Then $\epsilon_{i,j}\epsilon_{p,q}=
\epsilon_{p,q}\epsilon_{i,j}$.
\end{lemma}

\begin{proof}
Since all elements $\epsilon_{i,j}$ are obtained 
from $\theta_1$ via automorphisms, it is enough to show that 
$\theta_1$ commutes with all elements $\epsilon_{i,j}$ such that
$\{i,j\}\cap\{1,2\}=\varnothing$. 
Take any $v\in S_n$ such that 
$v(\{1,2\})=\{1,2\}$ and
$v(\{i,j\})=\{3,4\}$.
Such $v$ obviously exists. Then
$\theta_1$ commutes with $\epsilon_{i,j}$ if and only
if $v^{-1}\theta_1 v=\theta_1$ commutes with 
$v^{-1}\epsilon_{i,j}v=\theta_{3}$. 
The statement now follows from \eqref{it21}.
\end{proof}

\begin{lemma}\label{lemma204}
Suppose $\{i,j\}\cap\{p,q\}\neq\varnothing$. Then $\epsilon_{i,j}\epsilon_{p,q}=u\theta_1v$
for certain $u,v\in S_n$.
\end{lemma}
\begin{proof} If $\{i,j\}=\{p,q\}$ the statement is obvious as 
$\epsilon_{i,j}$ is an idempotent. Assume 
$\vert\{i,j\}\cap\{p,q\}\vert=1$. Since all elements 
$\epsilon_{i,j}$ are obtained  from $\theta_1$ via automorphisms, 
it is enough to consider the case when $\{i,j\}=\{1,2\}$, $p=2$
and $q>2$. Consider $v\in S_n$ such that  $v(1)=1, v(2)=2$ and $v(q)=3$.
Then, using  \eqref{it31}, \eqref{it11} and \eqref{it41} we have
\begin{displaymath}
v^{-1}\theta_1\epsilon_{p,q}v=\theta_1\theta_2=
\theta_1\sigma_1\theta_2\sigma_1\sigma_1=
\theta_1\sigma_2\theta_1\sigma_2\sigma_1=
\theta_1\sigma_2\sigma_1.
\end{displaymath}
The statement follows.
\end{proof}

For each $k$, $1\leq k\leq [\frac{n}{2}]$, set
$\delta_k=\theta_1\theta_3\dots\theta_{2k-1}$. Set also
$\delta_0=e$. The elements $\delta_i$, $0\leq i\leq [\frac{n}{2}]$, 
will be called {\em canonical}. The group $S_n\times S_n$ acts 
naturally on $T$  via $(g,h)(x)=g^{-1}xh$ for $x\in T$ and 
$(g,h)\in S_n\times S_n$. 

\begin{lemma}\label{lemma207}
Every $S_n\times S_n$-orbit contains a
canonical element.
\end{lemma}

\begin{proof}
Let $x\in T$. If $x\in S_n$ the statement is obvious. 
Assume that $x\not\in S_n$.
By Lemma \ref{lemma202} we can write $x=w\theta_1g_1\theta_1g_2\dots \theta_1g_k$ for some
$k\geq 1$ and $w, g_1,\dots, g_{k}\in S_n$. Moreover, we 
may assume that $x$ can not be written as a product of
$\theta_1$'s and elements of $S_n$, which contains less than 
$k$ occurrences of $\theta_1$. We have
\begin{multline}\label{eq:a3}
x=w(g_1\dots g_{k})(g_1\dots g_{k})^{-1}\theta_1(g_1\dots g_{k})\cdot
\\ \cdot
(g_2\dots g_{k})^{-1}\theta_1
(g_2\dots g_{k}) 
\dots (g_{k-1}g_{k})^{-1}\theta_1 (g_{k-1}g_{k}) g_{k}^{-1}\theta_1 g_{k},
\end{multline}
and hence we can write 
\begin{equation}\label{anya}
x=u\epsilon_{i_1,j_1}\dots \epsilon_{i_k,j_k},
\end{equation}
where $u=wg_1\dots g_{k}$ and $\{i_t,j_t\}$=$\{(g_t\dots g_k)(1),$ $(g_t\dots g_k)(2)\}$,
$1\leq t\leq k$. Since $x$ is chosen such that it can not be
reduced to an element of $T$ which contains less that $k$ entries of $\theta_1$, from Lemma~\ref{lemma203} and Lemma~\ref{lemma204}
it follows  that $\{i_t,j_t\}\cap\{i_s,j_s\}=\varnothing$ for 
any two factors  $\epsilon_{i_t,j_t}$, $\epsilon_{i_s,j_s}$ in 
\eqref{anya}. This implies that the $S_n\times S_n$-orbit of $x$ contains
$\epsilon_{i_1,j_1}\dots \epsilon_{i_k,j_k}$ with
$\{i_t,j_t\}\cap\{i_s,j_s\}=\varnothing$ for all $s\neq t$.

Now consider some $v\in S_n$ such that $v(i_1)=1$, $v(j_1)=2$,
$v(i_2)=3$ and so on, $v(j_k)=2k$. Then the element
$v^{-1}\epsilon_{i_1,j_1}\cdots \epsilon_{i_k,j_k}v$ is canonical
by definition. This completes the proof.
\end{proof}

\begin{remark}\label{anya1}
{\rm
From the proof of Lemma \ref{lemma207} it follows that each $x\in T$ 
can be written in the form 
$x=w\theta_1g_1\theta_1g_2\dots \theta_1g_k$, where 
$k\leq \lfloor\frac{n}{2}\rfloor$.
}
\end{remark}

\begin{lemma}\label{lemma222}
The $S_n\times S_n$-orbit of the canonical element 
$\delta_k$, $0\leq k\leq [\frac{n}{2}]$, contains at most 
\begin{displaymath}
\frac{(n!)^2}{2^{2k}(k!)^2(n-2k)!}
\end{displaymath}
elements.
\end{lemma}

\begin{proof}
It is enough to show that the stabilizer of $\delta_k$ under the $S_n\times S_n$-action
contains at least $(k!)^22^{2k}(n-2k)!$ elements.  
Set 
\begin{gather*}\Sigma^0_i=\sigma_{2i}\sigma_{2i-1}\sigma_{2i+1}
\sigma_{2i},\,\, 1\leq i\leq k-1;\\
\Sigma^1_i=\sigma_{2i}\sigma_{2i-1}\sigma_{2i+1}\sigma_{2i}
\sigma_{2i-1},\,\,\,\, 1\leq i\leq k-1.
\end{gather*}
Then both $\Sigma^0_i$  and $\Sigma^1_i$ swap 
the sets $\{2i-1, 2i\}$ and $\{2i+1, 2i+2\}$. It follows that 
the group $H$, generated by all $\Sigma^0_i$, consists of all 
permutations of the set $\{1,2\}, \{3,4\},\dots, \{2k-1, 2k\}$ 
and is therefore isomorphic to the group $S_k$. It is further
easy to see that the group $\tilde{H}$, generated by all $\Sigma^0_i$ and 
$\Sigma^1_i$, is isomorphic to the wreath product $H\wr S_2$.
From \eqref{it51} and \eqref{it31} it follows that the left
multiplication with both $\Sigma^0_i$ and $\Sigma^1_i$ 
stabilizes $\delta_k$. Therefore for each element of $\tilde{H}$
the left multiplication with this element stabilizes 
$\delta_k$ as well. Similarly one proves that the right 
multiplication with each element from $\tilde{H}$
stabilizes $\delta_k$.
Apart from this, from \eqref{it31} we have that
the conjugation by any element from the group 
$H'=\langle \sigma_{2k+1}, \dots, \sigma_{n-1}\rangle\simeq S_{n-2k}$
stabilizes $\delta_k$.

Observe that the group, generated by the left copy of
$\tilde{H}$, the right copy of $\tilde{H}$, and the $H'$ is a 
direct product of these three componets.  Using the product rule we 
derive that the cardinality of the stabilizer of $\delta_k$ is at least 
\begin{displaymath}
(\vert H\wr S_2\vert)^2\vert S_{n-2k}\vert= (k!)^22^{2k}(n-2k)!,
\end{displaymath}
and the proof is complete.
\end{proof}

\begin{corollary}\label{corollary31}
\begin{displaymath}
\vert T\vert \leq 
\sum_{k=0}^{\lfloor\frac{n}{2}\rfloor}\frac{(n!)^2}{2^{2k}(k!)^2(n-2k)!}.
\end{displaymath}
\end{corollary}

\begin{proof}
The proof follows from Lemma \ref{lemma222} and Remark~\ref{anya1}
by a direct calculation.
\end{proof}

\begin{proof}[Proof of Theorem~\ref{theorem1}]
Comparing Corollary~\ref{corollary31} and
Proposition~\ref{comb}\eqref{comb.2} we have
$|T|\leq |\mathfrak{B}_n|$. Since $\varphi:T\to \mathfrak{B}_n$
is surjective we have $|T|\geq |\mathfrak{B}_n|$. Hence
$|T|=|\mathfrak{B}_n|$ and $\varphi$ is an isomorphism.
\end{proof}

\section{Presentation for $\mathcal{IT}_n$}\label{s4}

For $i\in\{1,2,\dots,n-1\}$ let $\varrho_{i}$ denote the
element $\{i,i+1,i',(i+1)'\}\cup\bigcup_{j\neq i,i+1}\{j,j'\}\in
\mathcal{IT}_n$. By \cite[Proposition~9]{Mal2}, the elements
$\{\sigma_i\}$ and $\{\varrho_{i}\}$ generate $\mathcal{IT}_n$
(and even $\{\sigma_i\}$ and, say  $\varrho_{1}$, do).

Let $T$ denote the monoid with the identity element $e$, 
generated by the elements $\sigma_i$, $\tau_i$, 
$i=1,\dots,n-1$, subject to the following relations
(where $i,j\in\{1,2,\dots,n-1\}$):
\begin{gather}
\sigma_i^2=e; \quad
\sigma_i\sigma_j=\sigma_j\sigma_i,\,\, |i-j|>1;\quad
\sigma_i\sigma_j\sigma_i=\sigma_j\sigma_i\sigma_j, \,\, |i-j|=1;\label{it1}\\
\tau_i^2=\tau_i;\quad
\tau_i\tau_j=\tau_j\tau_i, \,\, i\neq j;\label{it2}\\
\tau_i\sigma_i=\sigma_i\tau_i=\tau_i; \quad
\tau_i\sigma_j=\sigma_j\tau_i, \,\, |i-j|>1;\label{it3}\\
\sigma_i\tau_j\sigma_i=\sigma_j\tau_i\sigma_j\,\,\text{ and }\,\,
\tau_i\sigma_j\tau_i=\tau_i\tau_j, \,\, |i-j|=1.\label{it4}
\end{gather}

\begin{theorem}\label{theorem2}
The map $\sigma_i\mapsto s_i$ and $\tau_i\to \varrho_i$, 
$i=1,\dots,n-1$,
extends to an isomorphism, $\varphi:T\to \mathcal{IT}_n$.
\end{theorem}

The rest of the section will be devoted to the proof of
Theorem~\ref{theorem2}.

It is a direct calculation to verify that the generators
$s_i$ and $\varrho_i$ of $\mathcal{IT}_n$ satisfy the relations,
corresponding to \eqref{it1}--\eqref{it4}. Thus the map
$\sigma_i\mapsto s_i$ and $\tau_i\mapsto \varrho_i$, $i=1,\dots,n-1$,
extends to an epimorphism, $\varphi:T\tto \mathcal{IT}_n$. Hence,
to prove Theorem~\ref{theorem2} we have only to show that 
$|T|=|\mathcal{IT}_n|$. As in the previous section, to do this we 
will study the structure of $T$ in details. 
Let $W$ denote the free monoid, generated by $\sigma_i$,
$\tau_i$, $i=1,\dots,n-1$,  $\psi:W\tto T$ denote the
canonical projection, and $\sim$ be the corresponding congruence on
$W$. The first part of our arguments is very similar to that from
the previous Section.

\begin{lemma}\label{lemma1}
The elements $\sigma_i$, $i=1,\dots,n-1$, generate the group 
$G$ of units in $T$, which is isomorphic to the symmetric group $S_n$
(and will be identified with $S_n$ in the sequel).
\end{lemma}

\begin{proof}
Analogous to that of Lemma~\ref{lemma301}.
\end{proof}

There are two natural actions on $T$:
\begin{enumerate}[(I)]
\item\label{actionone} The group $S_n$ acts on $T$ by inner 
automorphisms via conjugation.
\item\label{actiontwo} The group $S_n\times S_n$ acts on $T$ 
via $(g,h)(x)=g^{-1}xh$ for $x\in T$ and $(g,h)\in S_n\times S_n$.
\end{enumerate}

\begin{lemma}\label{lemma101}
The $S_n$-stabilizer of $\tau_1$ is the subgroup $H$ of
$S_n$, consisting of all permutations, which preserve the set
$\{1,2\}$. This subgroup is isomorphic to $S_2\times S_{n-2}$.
\end{lemma}

\begin{proof}
Analogous to that of Lemma~\ref{lemma201}.
\end{proof}

Since $S_n$ acts on $T$ via automorphisms and $\tau_1$ is an 
idempotent, all elements in the $S_n$-orbit of $\tau_1$ are idempotents.
From Lemma~\ref{lemma101} it follows that the elements of the
$S_n$-orbit of $\tau_1$ are in the natural bijection with the
cosets $H\backslash S_n$. By the definition of $H$, two elements,
$x,y\in S_n$, are contained in the same coset if and only if 
$x(\{1,2\})=y(\{1,2\})$.

\begin{lemma}\label{lemma102}
The $S_n$-orbit of $\tau_1$ contains all $\tau_i$,
$i=1,\dots,n-1$. Moreover, for $w\in S_n$ we have
$w^{-1}\tau_1 w=\tau_i$ if and only if $w(\{1,2\})=\{i,i+1\}$.
\end{lemma}

\begin{proof}
Analogous to that of Lemma~\ref{lemma202}.
\end{proof}

\begin{lemma}\label{lemma103}
All elements in the $S_n$-orbit of $\tau_1$ commute.
\end{lemma}

\begin{proof}
Since all elements in the $S_n$-orbit of $\tau_1$ are obtained 
from $\tau_1$ via automorphisms, it is enough to show that 
$\tau_1$ commutes with all elements in this orbit. 
Let $w\in S_n$ be such that $w(\{1,2\})=\{i,j\}$.
If $\{i,j\}=\{1,2\}$ then $w^{-1}\tau_1 w=\tau_1$
by Lemma~\ref{lemma102} and hence we may assume $\{i,j\}\neq \{1,2\}$. 

Take any $v\in S_n$ such that 
\begin{itemize}
\item $v(\{1,2\})=\{1,2\}$ and
$v(\{i,j\})=\{3,4\}$ if $\{i,j\}\cap\{1,2\}=\varnothing$;
\item $v(\{1,2\})=\{1,2\}$ and
$v(\{i,j\})=\{2,3\}$ if $\{i,j\}\cap\{1,2\}\neq \varnothing$.
\end{itemize}
Such $v$ obviously exists. Then
$\tau_1$ commutes with $w^{-1}\tau_1 w$ if and only
if $v^{-1}\tau_1 v$ commutes with 
$v^{-1}w^{-1}\tau_1 wv$. Using our choice of $v$ and
Lemma~\ref{lemma102} we have
$v^{-1}\tau_1 v=\tau_1$ and
$v^{-1}w^{-1}\tau_1 wv=\tau_j$, where $j=3$ if
$\{i,j\}\cap\{1,2\}=\varnothing$, and $j=2$ otherwise.
The statement now follows from \eqref{it2}.
\end{proof}

For $w\in S_n$ such that $w(\{1,2\})=\{i,j\}$, where $i<j$, we
set $\varepsilon_{i,j}=w^{-1}\tau_1 w$, which is well defined
by Lemma~\ref{lemma101}.

\begin{lemma}\label{lemma104}
Let $\{i,j,k\}\subset \{1,2,\dots,n\}$ and 
$i<j<k$. Then  
\begin{displaymath}
\varepsilon_{i,j}\varepsilon_{j,k}=
\varepsilon_{i,k}\varepsilon_{j,k}=
\varepsilon_{i,j}\varepsilon_{i,k}.
\end{displaymath}
\end{lemma}

\begin{proof}
We prove that $\varepsilon_{i,j}\varepsilon_{j,k}=
\varepsilon_{i,k}\varepsilon_{j,k}$ and the second equality is proved by
analogous arguments. Let  $w\in S_n$ be such that $w(i)=1$, $w(j)=2$, 
$w(k)=3$. Conjugating by $w$ we reduce our equality to the equality
$\tau_1\tau_2=\sigma_2\tau_1\sigma_2\tau_2$. 
Using \eqref{it4} twice and \eqref{it3} we have
\begin{displaymath}
\sigma_2\tau_1\sigma_2\tau_2=
\sigma_1\tau_2\sigma_1\tau_2=
\sigma_1\tau_1\tau_2=
\tau_1\tau_2.
\end{displaymath}
The claim follows.
\end{proof}

For $i,j\in M$ set $\varepsilon_{i,i}=e$ and
$\varepsilon_{i,j}=\varepsilon_{j,i}$ if $j<i$.
For a non-empty binary relation, $\rho$, on $M$ set
\begin{displaymath}
\varepsilon_{\rho}=\prod_{i\rho j}\varepsilon_{i,j}.
\end{displaymath}

\begin{corollary}\label{corollary71}
Let $\rho$ be non-empty binary relation on $M$ and
$\rho^*$ be the reflexive-symmetric-transitive closure of $\rho$.
Then $\varepsilon_{\rho}=\varepsilon_{\rho^*}$
\end{corollary}

\begin{proof}
Follows easily from Lemma~\ref{lemma103},
Lemma~\ref{lemma104} and the fact that all 
$\varepsilon_{i,j}$'s are idempotents.
\end{proof}

Let $\lambda:\{1,\dots,n\}=X_1\cup \dots \cup X_k$ be a 
decomposition of $M$ into an unordered union of pairwise
disjoint sets. With this decomposition we associate
the equivalence relation $\rho_{\lambda}$ on $M$, whose equivalence
classes coincide with $X_i$'s. 

\begin{corollary}\label{corollary72}
Let $\lambda$ and $\mu$ be two decompositions of $M$ as above.
Assume that the types of $\lambda$ and $\mu$
coincide. Then $\varepsilon_{\rho_\lambda}$ and 
$\varepsilon_{\rho_\mu}$ are conjugate in $T$.
\end{corollary}

\begin{proof}
Let $v\in S_n$ be an element, which maps $\lambda$ to $\mu$ (such element
exists since the types of $\lambda$ and $\mu$ are the same). One 
easily sees that $v^{-1}\varepsilon_{\rho_\lambda}v=\varepsilon_{\rho_\mu}$.
The statement follows.
\end{proof}

A decomposition, $\lambda:\{1,\dots,n\}=X_1\cup \dots \cup X_k$,
is called {\em canonical} provided that (up to a permutation of the blocks)
we have $|X_1|\geq |X_2|\geq \dots\geq |X_k|$,
$X_1=\{1,2,\dots,l_1\}$, $X_2=\{l_1+1,l_1+2,\dots,l_1+l_2\}$ and so on.
Note that in this case $\lambda$ can also be viewed as a  
{\em partition} of $n$. The element $\varepsilon_{\rho_\lambda}$
will be called {\em canonical} provided that $\lambda$ is canonical.

\begin{lemma}\label{lemma107}
Every $S_n\times S_n$-orbit contains a
canonical element.
\end{lemma}

\begin{proof}
Because of Corollary~\ref{corollary72} it is enough to show that 
every $S_n\times S_n$-orbit contains $\varepsilon_{\rho_\lambda}$
for some decomposition $\lambda$. Let $x\in T$. If
$x\in S_n$, then the statement is obvious. Let $x\in T\setminus S_n$.
From Lemma~\ref{lemma102} we have that the semigroup $T$ is generated
by $S_n$ and $\tau_1$. Hence we have
$x=w\tau_1g_1\tau_1g_2\cdots\tau_1g_k$ for some 
$w,g_1,\dots,g_k\in S_n$. Therefore
\begin{multline*}
x=w(g_1\dots g_{k})(g_1\dots g_{k})^{-1}\tau_1(g_1\dots g_{k})\cdot
\\ \cdot
(g_2\dots g_{k})^{-1}\tau_1
(g_2\dots g_{k}) 
\dots (g_{k-1}g_{k})^{-1}\tau_1 (g_{k-1}g_{k}) g_{k}^{-1}\tau_1 g_{k},
\end{multline*}
and hence we can write $x=u\varepsilon_{i_1,j_1}\dots \varepsilon_{i_k,j_k}$,
where $u=wg_1\dots g_{k}$ and 
\begin{displaymath}
\{i_t,j_t\}=\{(g_t\dots g_k)(1), 
(g_t\dots g_k)(2)\},\,\, 1\leq t\leq k.
\end{displaymath}
Define the equivalence relation $\rho$ as 
the reflexive-symmetric-transitive closure of the relation
$\{(i_1,j_1),\dots,(i_k,j_k)\}$ and let $\lambda$ be the corresponding
decomposition of $\{1,2,\dots,n\}$. From Corollary~\ref{corollary71}
we get that the  $S_n\times S_n$-orbit of $x$ contains
$\varepsilon_{\rho}=\varepsilon_{\rho_\lambda}$. This
completes the proof.
\end{proof}

\begin{lemma}\label{lemma108}
Let $\lambda$ be a canonical decomposition of 
$\{1,2,\dots,n\}$. For $i=1,\dots,n$ set
$\lambda^{(i)}=|\{j:|X_j|=i\}|$.
Then the $S_n\times S_n$-stabilizer of 
$\varepsilon_{\rho_{\lambda}}$ contains at least
\begin{displaymath}
\prod_{i=1}^n(\lambda^{(i)}!(i!)^{2\lambda^{(i)}})
\end{displaymath}
elements.
\end{lemma}

\begin{proof}
Fix $i\in\{1,2,\dots,n\}$.
Let  $X_{a},X_{a+1}\dots,X_{b}$ be all blocks of $\lambda$ of
cardinality $i$. Then for any non-maximal element $j$ of any
of $X_{a},X_{a+1}\dots,X_{b}$, using Lemma~\ref{lemma103}, the
definition of $\varepsilon_{\rho_{\lambda}}$, and \eqref{it3} we have
$\sigma_j\varepsilon_{\rho_{\lambda}}=
\varepsilon_{\rho_{\lambda}}\sigma_j=
\varepsilon_{\rho_{\lambda}}$. Moreover, for any 
$w\in S_n$, which stabilizes all elements outside
$X_{a}\cup X_{a+1}\cup \cdots\cup X_{b}$ and maps
each $X_s$ to some $X_t$, we have $w(\lambda)=\lambda$ and
hence $w^{-1}\varepsilon_{\rho_{\lambda}}w=\varepsilon_{\rho_{\lambda}}$.
This gives us exactly $\lambda^{(i)}!(i!)^{2\lambda^{(i)}}$
elements of the $S_n\times S_n$-stabilizer. The statement of the lemma
now follows by applying the product rule since for different $i$
the nontrivial elements $w$ above stabilize pairwise different 
subsets of $\{1,\dots,n\}$.
\end{proof}

\begin{corollary}\label{corollary9}
\begin{displaymath}
|T|\leq
\sum_{\lambda\vdash n}
\frac{(n!)^2}
{\displaystyle\prod_{i=1}^n(\lambda^{(i)}!(i!)^{2\lambda^{(i)}})}.
\end{displaymath}
\end{corollary}

\begin{proof}
Canonical elements of $T$ are in bijection with partitions
$\lambda\vdash n$ by construction. By Lemma~\ref{lemma107}, every
$S_n\times S_n$-orbit contains a canonical element.
We have $|S_n\times S_n|=(n!)^2$. By Lemma~\ref{lemma108}, the
stabilizer of a canonical element, corresponding to $\lambda$,
contains at least $\prod_{i=1}^n(\lambda^{(i)}!(i!)^{2\lambda^{(i)}})$
elements. The statement now follows by applying the sum rule.
\end{proof}

\begin{proof}[Proof of Theorem~\ref{theorem2}.]
Comparing Corollary~\ref{corollary9} and
Proposition~\ref{comb}\eqref{comb.3} we have
$|T|\leq |\mathcal{IT}_n|$. Since $\varphi:T\to \mathcal{IT}_n$
is surjective we have $|T|\geq |\mathcal{IT}_n|$. Hence
$|T|=|\mathcal{IT}_n|$ and $\varphi$ is an isomorphism.
\end{proof}

\begin{remark}\label{remark1}
{\rm
From the above arguments it follows that the inequality obtained
in Lemma~\ref{lemma108} is in fact an equality. From the
proof of Lemma~\ref{lemma108} one easily derives that the
$S_n\times S_n$-stabilizer of $\varepsilon_{\rho_{\lambda}}$ is
isomorphic to the direct product of wreath products
$S_{\lambda^{(i)}}\wr (S_i\times S_i)$. 
}
\end{remark}

\begin{remark}\label{remark2}
{\rm
Following the arguments of the proof of Theorem~\ref{theorem2}
one easily proves the following presentation for the
symmetric inverse semigroup $\mathcal{IS}_n$:
$\mathcal{IS}_n$ is generated, as a monoid, by
$\sigma_1,\dots,\sigma_{n-1},\vartheta_1,\dots,\vartheta_n$
subject to the following relations:
\begin{gather}
\sigma_i^2=e; \quad
\sigma_i\sigma_j=\sigma_j\sigma_i,\,\, |i-j|>1;\quad
\sigma_i\sigma_j\sigma_i=\sigma_j\sigma_i\sigma_j, \,\, |i-j|=1; 
\label{it50}\\
\vartheta_i^2=\vartheta_i;\quad
\vartheta_i\vartheta_j=\vartheta_j\vartheta_i\,\, i\neq j; 
\label{it42}\\
\sigma_i\vartheta_i=\vartheta_{i+1}\sigma_i; \,\,
\sigma_i\vartheta_j=\vartheta_{j}\sigma_i, \,\, j\neq i,i+1;\quad
\vartheta_i\sigma_i\vartheta_i=\vartheta_i\vartheta_{i+1}. 
\label{it43}
\end{gather}
The classical presentation for $\mathcal{IS}_n$ usually involves only
one additional generator (namely $\vartheta_1$) and can be 
found for example in \cite[Chapter~9]{Li}.
}
\end{remark}

\section{Presentation for $\mathcal{P}\mathfrak{B}_n$}\label{s5}

For $i\in\{1,\dots,n\}$ let
$\varsigma_i$ denote the element
$\{i\}\cup \{i'\}\cup \bigcup_{j\neq i}\{j,j'\}$. Using
\cite{Maz0}, it is easy to see that $\mathcal{P}\mathfrak{B}_n$
is generated by $\{\sigma_i\}\cup\{\pi_i\}\cup\{\varsigma_i\}$
(and even by $\{\sigma_i\}$, $\pi_1$ and $\varsigma_1$). 

Let $T$ denote the monoid with the identity element $e$, 
generated by the elements $\sigma_i$, $\theta_i$, 
$i=1,\dots,n-1$, and $\vartheta_i$, $i=1,\dots,n$, 
subject to the relations
\eqref{it11}--\eqref{it41new}, the relations from
Remark~\ref{remark2}, and the following relations
(for all appropriate $i$ and $j$):
\begin{gather}
\theta_i\vartheta_j=\vartheta_j\theta_i,\,\, j\neq i,i+1;\label{pb1}\\
\theta_i\vartheta_i=\theta_i\vartheta_{i+1}=
\theta_i\vartheta_i\vartheta_{i+1},\quad
\vartheta_i\theta_i=\vartheta_{i+1}\theta_i=
\vartheta_i\vartheta_{i+1}\theta_i;\label{pb2}\\
\theta_i\vartheta_i\theta_i=\theta_i,\quad
\vartheta_i\theta_i\vartheta_i=\vartheta_i\vartheta_{i+1}.\label{pb3}
\end{gather}

\begin{theorem}\label{theorem3}
The map $\sigma_i\mapsto s_i$, $\theta_i\to \pi_i$, $i=1,\dots,n-1$,
and $\vartheta_i\mapsto \varsigma_i$, $i=1,\dots,n$,
extends to an isomorphism, $\varphi:T\to \mathcal{P}\mathfrak{B}_n$.
\end{theorem}

We will again start with the following auxiliary technical statement,
which we will need later:

\begin{lemma}\label{technicallemma}
Under the assumption that \eqref{it11}--\eqref{it41new}, \eqref{pb1}--\eqref{pb3} and the relations from
Remark~\ref{remark2} are satisfied, one has the relation
\begin{equation}
\sigma_{i+2}\sigma_{i+1}\theta_i\vartheta_{i+2}\vartheta_{i+3}=
\sigma_{i}\sigma_{i+1}\vartheta_i\theta_i\theta_{i+2}\vartheta_{i+2}.
\label{pb4}
\end{equation}
\end{lemma}

\begin{proof}
Using \eqref{it41new} twice and \eqref{it11} we have
\begin{multline*}
\sigma_{i+2}\sigma_{i+1}\theta_i\vartheta_{i+2}\vartheta_{i+3}=
\sigma_{i+2}\sigma_{i}\theta_{i+1}\theta_i\vartheta_{i+2}\vartheta_{i+3}=\\=
\sigma_{i}\sigma_{i+2}\theta_{i+1}\theta_i\vartheta_{i+2}\vartheta_{i+3}=
\sigma_{i}\sigma_{i+1}\theta_{i+2}\theta_{i+1}
\theta_i\vartheta_{i+2}\vartheta_{i+3}, 
\end{multline*}
and hence \eqref{pb4} reduces to 
\begin{equation}
\theta_{i+2}\theta_{i+1}\theta_i\vartheta_{i+2}\vartheta_{i+3}=
\vartheta_i\theta_i\theta_{i+2}\vartheta_{i+2}.
\label{pb4b}
\end{equation}
Using \eqref{pb1}-\eqref{pb3} and \eqref{it21} we have
\begin{multline*}
\theta_{i+2}\theta_{i+1}\theta_i\vartheta_{i+2}\vartheta_{i+3}=
\theta_{i+2}\vartheta_{i+3}\theta_{i+1}\vartheta_{i+2}\theta_i=
\theta_{i+2}\vartheta_{i+2}\theta_{i+1}\vartheta_{i+1}\theta_i=\\=
\theta_{i+2}\vartheta_{i+1}\theta_{i+1}\vartheta_{i+1}\theta_i=
\theta_{i+2}\vartheta_{i+2}\vartheta_{i+1}\theta_i=
\vartheta_i\theta_i\theta_{i+2}\vartheta_{i+2},
\end{multline*}
which gives \eqref{pb4b}. The statement follows.
\end{proof}

As in the previous section, one easily checks that this map extends
to an epimorphism and hence to complete the proof one has to compare
the cardinalities of $T$ and $\mathcal{P}\mathfrak{B}_n$.

Similarly to what was done in Section \ref{s4}, using the presentation 
of $\mathcal{IS}_n$ given in Remark \ref{remark2}, one proves that 
elements $\sigma_i$, $i=1,\dots,n-1$, generate  the symmetric group 
$S_n$, and that the elements $\sigma_i$,  $i=1,\dots,n-1$; 
$\vartheta_i$, $i=1,\dots, n$, generate the semigroup, which is 
isomorphic to $\mathcal{IS}_n$ (and which will be identified with it).
As in  Section \ref{s4}  we consider the natural action of $S_n$ on 
$T$ by inner automorphisms of $T$ via conjugation: $x^g=g^{-1}xg$ for 
each $x\in T$, $g\in S_n$.  Set $\xi_i=\theta_i\vartheta_i$,
$\eta_i=\vartheta_i\theta_i$, $1\leq i\leq n-1$.

\begin{lemma}\label{lemma401}
The $S_n$-stabilizer of each of $\theta_1$, 
$\xi_1$, $\eta_1$ is the subgroup $H$ of
$S_n$, consisting of all permutations, which preserve the set
$\{1,2\}$. This subgroup is isomorphic to $S_2\times S_{n-2}$.
\end{lemma}

\begin{proof}
For $\theta_1$ this follows from Lemma \ref{lemma201}.
For each $j\geq 2$ we have that $\sigma_j$ commutes with both 
$\xi_1$ and $\eta_1$ by \eqref{it31} and \eqref{it43} respectively, 
and hence $\sigma_j\xi_1\sigma_j=\xi_1$ and  
$\sigma_j\eta_1\sigma_j=\eta_1$.
Let $j=1$. Then
\begin{gather*}
\sigma_1\xi_1\sigma_1=\sigma_1\theta_1\vartheta_1\sigma_1=
\sigma_1\theta_1\sigma_1\vartheta_2=\theta_1\vartheta_2=
\theta_1\vartheta_1=\xi_1;\\
\sigma_1\eta_1\sigma_1=\sigma_1\vartheta_1\theta_1\sigma_1=
\vartheta_2\sigma_1\theta_1\sigma_1=\vartheta_2\theta_1=
\vartheta_1\theta_1=\eta_1
\end{gather*}
by \eqref{it43} and \eqref{it31}. Hence $\sigma_1$ also stabilizes 
$\xi_1$ and $\eta_1$. Since $\sigma_j$,  $j\neq 2$, generate $H$, we 
obtain that all elements of $H$ stabilize $\xi_1$ and $\eta_1$. 
In particular, the $S_n$-orbits of $\xi_1$ and of $\eta_1$ consist 
of at most $|S_n|/|H|=\binom{n}{2}$ elements each. At the same time, 
the $S_n$-orbits of $\varphi(\xi_1)$ and  $\varphi(\eta_1)$ consist 
of exactly $\binom{n}{2}$ different elements and hence $H$ must 
coincide with the $S_n$-stabilizer of both $\xi_1$ and $\eta_1$.
\end{proof}

Since $S_n$ acts on $T$ via automorphisms and $\theta_1$, $\xi_1$, $\eta_1$ are idempotents,
all elements in the $S_n$-orbits of $\theta_1$, $\xi_1$, $\eta_1$ are idempotents as well.
From Lemma~\ref{lemma401} it follows that the elements of the
$S_n$-orbits of $\theta_1$, $\xi_1$, $\eta_1$ are in the natural bijections with the
cosets $H\backslash S_n$. By the definition of $H$, two elements,
$x,y\in S_n$, are contained in the same coset if and only if 
$x(\{1,2\})=y(\{1,2\})$.

\begin{lemma}\label{lemma402}
The $S_n$-orbits of $\theta_1$, $\xi_1$, $\eta_1$ contain all 
elements $\theta_i$, $\xi_i$ and $\eta_i$,  $i=1,\dots, n-1$, 
respectively. Moreover, for $w\in S_n$ we have 
$w^{-1}\theta_1 w=\theta_i$ if and only if $w(\{1,2\})=\{i,i+1\}$ and
analogously for $\xi_1$ and $\eta_1$.
\end{lemma}

\begin{proof}
The proof for the $S_n$-orbit of $\theta_1$ is analogous to that of 
Lemma \ref{lemma202}. We prove the statement for the $S_n$-orbit of 
$\xi_1$. For the $S_n$-orbit of $\eta_1$ the arguments are analogous. 
We use induction on $i$ with the case $i=1$ being trivial.
Let $i>1$ and assume that $\xi_{i-1}$ is contained in our orbit.
Then, using \eqref{it43}, \eqref{it11} and \eqref{it41}, we compute
\begin{multline*}
\xi_{i}=\theta_i\vartheta_i=
\sigma_{i-1}\sigma_{i}\theta_{i-1}\sigma_{i}\sigma_{i-1}\vartheta_i=
\sigma_{i-1}\sigma_{i}\theta_{i-1}\sigma_{i}\vartheta_{i-1}\sigma_{i-1}= \\ \sigma_{i-1}\sigma_{i}
\theta_{i-1}\vartheta_{i-1}\sigma_{i}\sigma_{i-1}=\sigma_{i-1}\sigma_{i}
\xi_{i-1}\sigma_{i}\sigma_{i-1},
\end{multline*}
and hence $\xi_{i}$ is contained in our orbit as well. 
The second claim follows from \eqref{oldequation}.
This completes the proof.
\end{proof}

For $w\in S_n$ such that $w(\{1,2\})=\{i,j\}$, where $i<j$, we
set $\epsilon_{i,j}=w^{-1}\theta_1 w$, 
$\mu_{i,j}=w^{-1}\xi_1 w$, $\nu_{i,j}=
w^{-1}\eta_1 w$. All these elements are well defined
by Lemma~\ref{lemma401}.

\begin{lemma}\label{lemma520}
\begin{enumerate}[(a)]
\item\label{lemma520.1} $\vartheta_i\epsilon_{i,j}=
\vartheta_j\epsilon_{i,j}=\vartheta_i\vartheta_j\epsilon_{i,j}=\nu_{i,j}$;
$\vartheta_k\epsilon_{i,j}=\epsilon_{i,j}\vartheta_k$, $k\not\in \{i,j\}$.
\item\label{lemma520.2} $\vartheta_i\mu_{i,j}=\vartheta_j\mu_{i,j}=\vartheta_i\vartheta_j\mu_{i,j}
=\vartheta_i\vartheta_j$;
$\vartheta_k\mu_{i,j}=\mu_{i,j}\vartheta_k$, $k\not\in \{i,j\}$.
\end{enumerate}
\end{lemma}

\begin{proof}
First we prove \eqref{lemma520.1}.
Because of Lemma~\ref{lemma402} it is enough to check that $\vartheta_1\epsilon_{1,2}=\vartheta_2\epsilon_{1,2}=
\vartheta_1\vartheta_2\epsilon_{1,2}=\nu_{1,2}$ and that
$\vartheta_3\epsilon_{1,2}=\epsilon_{1,2}\vartheta_3$. The latter 
equalities follow from \eqref{pb2} and \eqref{pb1}.

Now we prove \eqref{lemma520.2}.
Again, because of Lemma~\ref{lemma402} it is enough to 
check that $\vartheta_1\mu_{1,2}=\vartheta_2\mu_{1,2}=
\vartheta_1\vartheta_2\mu_{1,2}=\vartheta_1\vartheta_2$ and that
$\vartheta_3\mu_{1,2}=\mu_{1,2}\vartheta_3$. Using \eqref{pb3}, 
\eqref{pb2} and \eqref{pb1} we have
\begin{displaymath}
\vartheta_1\mu_{1,2}=\vartheta_1\theta_1\vartheta_1=
\vartheta_1\vartheta_2; \quad
\vartheta_1\mu_{2,3}=\vartheta_1\theta_2\vartheta_2=
\theta_2\vartheta_1\vartheta_2=\theta_2\vartheta_2
\vartheta_1=\mu_{2,3}\vartheta_1,
\end{displaymath}
as required.
\end{proof}

\begin{lemma}\label{lemma403}
Suppose $\{i,j\}\cap\{p,q\}=\varnothing$. Then $\epsilon_{i,j}\epsilon_{p,q}=
\epsilon_{p,q}\epsilon_{i,j}$, $\mu_{i,j}\mu_{p,q}=
\mu_{p,q}\mu_{i,j}$ and $\epsilon_{i,j}\mu_{p,q}=
\mu_{p,q}\epsilon_{i,j}$.
\end{lemma}

\begin{proof}
Following the arguments from the proof of Lemma \ref{lemma203} it 
is enough to show that $\mu_{1,2}\mu_{3,4}=
\mu_{3,4}\mu_{1,2}$ and $\mu_{1,2}\epsilon_{3,4}=
\epsilon_{3,4}\mu_{1,2}$, that is that $\xi_1\xi_3=\xi_3\xi_1$ 
and $\xi_1\theta_3=\theta_3\xi_1$.
Using \eqref{pb1}, \eqref{it42} and \eqref{it21} we have
\begin{displaymath}
\xi_1\xi_3=\theta_1\vartheta_1\theta_3\vartheta_3=
\theta_1\theta_3\vartheta_1\vartheta_3=
\theta_3\theta_1\vartheta_3\vartheta_1=
\theta_3\vartheta_3\theta_1\vartheta_1=\xi_3\xi_1,
\end{displaymath}
and using \eqref{pb1} and \eqref{it21} we also obtain
$\xi_1\theta_3=\theta_1\vartheta_1\theta_3=
\theta_1\theta_3\vartheta_1=\theta_3\xi_1$,
as required.
\end{proof}

\begin{lemma}\label{lemma404}
Suppose $\{i,j\}\cap\{p,q\}\neq\varnothing$. Then 
each of the elements $\epsilon_{i,j}\epsilon_{p,q}$,
$\mu_{i,j}\mu_{p,q}$, $\epsilon_{i,j}\mu_{p,q}$, $\mu_{i,j}\epsilon_{p,q}$
equals to the element of the form $u\theta_1v$ for some
$u,v\in \mathcal{IS}_n$. 
\end{lemma}

\begin{proof} 
Using the argument from the proof of Lemma \ref{lemma204} it is enough to 
prove the statement only for the elements $\mu_{1,2}\mu_{2,3}$,
$\mu_{1,2}\epsilon_{2,3}$, $\epsilon_{1,2}\mu_{2,3}$. We have 
\begin{displaymath}
\mu_{1,2}\mu_{2,3}=\xi_1\xi_2=\theta_1\vartheta_1\theta_2\vartheta_2=
\theta_1\vartheta_2\theta_2\vartheta_2=\theta_1\vartheta_2\vartheta_3=
\xi_1\vartheta_3=\theta_1\vartheta_1\vartheta_3
\end{displaymath}
by \eqref{pb2} and \eqref{pb3}; and
\begin{multline*}
\mu_{1,2}\epsilon_{2,3}=\theta_1\vartheta_1\theta_2=
\theta_1\vartheta_1\sigma_1\sigma_2\theta_1
\sigma_1\sigma_2=\theta_1\sigma_1\vartheta_2\sigma_2\theta_1
\sigma_1\sigma_2=\\
\theta_1\sigma_1\sigma_2\vartheta_3\theta_1
\sigma_1\sigma_2=\theta_1\sigma_1\sigma_2\theta_1\vartheta_3
\sigma_1\sigma_2=\theta_1
\sigma_2\theta_1\vartheta_3\sigma_1\sigma_2=
\theta_1\vartheta_3\sigma_1\sigma_2
\end{multline*}
by \eqref{it11}, \eqref{it41}, \eqref{it31}, \eqref{it43}.
Finally, 
\begin{displaymath}
\epsilon_{1,2}\mu_{2,3}=\theta_1\theta_2\vartheta_2=
\theta_1\sigma_1\sigma_2\theta_1\sigma_2\sigma_1\vartheta_2=
\theta_1\sigma_2\vartheta_1\sigma_1.
\end{displaymath}
using \eqref{it11}, \eqref{it31} and \eqref{it41}.
The statement follows.
\end{proof}

For each subset $\{i_1,\dots, i_k\}$ of $\{1,2,\dots, n\}$ 
set $\vartheta(\{i_1,\dots, i_k\})=
\vartheta_{i_1}\dots \vartheta_{i_k}$.
Obviously, $\vartheta(\{i_1,\dots, i_k\})$ is an idempotent and
each idempotent of $\mathcal{IS}_n$ has such a form.
In the sequel we will use the obvious fact that each element 
of $\mathcal{IS}_n$ can be written in the form $uv$, where
$u$ is an idempotent, and $v\in S_n$. 

As in the previous sections we consider the
$S_n\times S_n$-action on $T$ given by 
$(g,h)(x)=g^{-1}xh$ for $x\in T$ and 
$(g,h)\in S_n\times S_n$. 

\begin{lemma}\label{lemma510}
Every $S_n\times S_n$-orbit contains either $e$ or an element of the form
$\vartheta(A)\gamma_{i_1,j_1}\dots\gamma_{i_s,j_s}$,
where $A\subset\{1,2,\dots, n\}$,
the sets $\{i_l,j_l\}$ are pairwise disjoint, and each $\gamma_{i_l,j_l}$
equals either $\epsilon_{i_l,j_l}$ or $\mu_{i_l,j_l}$.
\end{lemma}

\begin{proof}
The idea of the proof is analogous to that of Lemma \ref{lemma207}.
Let $x\in T$. If $x\in S_n$ the statement is obvious. 
Assume that $x\not\in S_n$. Since $T$ is generated by 
${\mathcal{IS}}_n$ and $\theta_1$ we can write 
\begin{equation}\label{eq:a}
x=wu\theta_1u_1g_1\theta_1u_2g_2\cdots \theta_1u_kg_k
\end{equation} 
for some $k\geq 1$, $w, g_1,\dots, g_{k}\in S_n$ and
$u, u_1,\dots, u_k\in E(\mathcal{IS}_n)$. Moreover, we 
may assume that $x$ can not be written as a product of
$\theta_1$'s and elements of ${\mathcal{IS}}_n$, which contains 
less than  $k$ occurrences of $\theta_1$. 
We claim that $x$ can be written as
\begin{equation}\label{eq:a1}
x=wu'\gamma^1_1g_1'\gamma^2_1g_2'\cdots \gamma^k_1g_k',
\end{equation}
where, $w, g_1',\dots, g_k'\in S_n$, $u'\in E({\mathcal{IS}}_n)$, 
and each $\gamma^i_1$ is equal to 
either $\theta_1$ or $\xi_1$. Let us prove this by
induction on $k$. Let $k=1$ and $x=wu\theta_1u_1g_1$. We know that
$u_1=\vartheta(B)$ for some $B\subset\{1,\dots, n\}$.
Let $A=B\setminus \{1,2\}$.
Using \eqref{pb1} and \eqref{pb2} we obtain
that
\begin{equation*}
x=\left\lbrace\begin{array}{l} wuu_1\theta_1 g_1, \text{ if }B\cap\{1,2\}=\varnothing;\\
wu\vartheta(A)\xi_1g_1, \text{ if }B\cap\{1,2\}\neq\varnothing,
\end{array}\right.
\end{equation*}
as required. Let now $k\geq 2$. Applying the basis of the induction to $\theta_1u_kg_k$
we obtain
\begin{multline*}
x=wu\theta_1u_1g_1\theta_1u_2g_2\cdots \theta_1u_{k-1}g_{k-1}\theta_1u_kg_k=\\
wu\theta_1u_1g_1\theta_1u_2g_2\cdots \theta_1u_{k-1}g_{k-1}u_k'\gamma^k_1g_k,
\end{multline*}
where $u_k'$ is an idempotent of ${\mathcal{IS}}_n$ and $\gamma^k_1$ is 
either $\xi_1$ or $\theta_1$. Now, since $u_{k-1}g_{k-1}u_k'\in
{\mathcal{IS}}_n$, we can write $u_{k-1}g_{k-1}u_k'=u_{k-1}'g_{k-1}'$ 
for some $g_{k-1}'\in S_n$ and $u_{k-1}'\in E({\mathcal{IS}}_n)$. 
Now \eqref{eq:a1} follows by applying the inductive assumption to
$wu\theta_1u_1g_1\theta_1u_2g_2\cdots u_{k-2}g_{k-2}
\theta_1u_{k-1}'g_{k-1}'$.

Similarly to \eqref{eq:a3} we can rewrite \eqref{eq:a1} as follows:
\begin{multline*}
x=wu'(g_1'\cdots g_{k}')(g_1'\cdots g_{k}')^{-1}
\gamma^1_1(g_1'\cdots g_{k}')\cdot
\\ \cdot
(g_2'\cdots g_{k}')^{-1}\gamma^2_1
(g_2'\cdots g_{k}') 
\cdots (g_{k-1}'g_{k}')^{-1}\gamma^{k-1}_1 (g_{k-1}'g_{k}') g_{k}'^{-1}\gamma^k_1 g_{k}',
\end{multline*}
and therefore we can write 
\begin{equation}\label{eq:a4}
x=vu'\gamma_{i_1,j_1}\cdots\gamma_{i_k,j_k},
\end{equation}
where  $v=wg_1'\cdots g_{k}'$, $\{i_t,j_t\}$=$\{(g_t'\cdots g_k')(1), 
(g_t'\cdots g_k')(2)\}$,
$1\leq t\leq k$, and each $\gamma_{i_l,j_l}$
is equal to either $\epsilon_{i_l,j_l}$ or $\mu_{i_l,j_l}$.
Since $x$ is initially chosen such that it can not be
reduced to an element of $T$, which contains less that $k$ entries of $\theta_1$, from Lemma~\ref{lemma404} 
it follows  that $\{i_t,j_t\}\cap\{i_l,j_l\}=\varnothing$ for 
any two factors $\gamma_{i_t,j_t}$, $\gamma_{i_l,j_l}$ in 
\eqref{eq:a4}. This implies that the $S_n\times S_n$-orbit of $x$ contains
$u'\gamma_{i_1,j_1}\cdots\gamma_{i_s,j_s}$ such that $u'\in E(\IS_n)$, 
$\{i_t,j_t\}\cap\{i_l,j_l\}=\varnothing$ for all $l\neq t$.
The statement follows.
\end{proof}

\begin{corollary}\label{cor530}
Any $S_n\times S_n$- orbit contains either $e$ or an element of the form
$\vartheta(A)\gamma_{i_1,j_1}\cdots\gamma_{i_s,j_s}$, such that
\begin{enumerate}[(i)]
\item the sets $\{i_l,j_l\}$ are pairwise disjoint;
\item each $\gamma_{i_l,j_l}$ equals to either $\epsilon_{i_l,j_l}$ or $\mu_{i_l,j_l}$ or $\nu_{i_l,j_l}$;
\item $A\cap\{i_1, j_1, \dots i_s, j_s\}=\varnothing$.
\end{enumerate}
\end{corollary}

\begin{proof}
This follows from Lemma \ref{lemma510} and Lemma \ref{lemma520}.
\end{proof}

Now we introduce the notion of a canonical element.
Let $k,l,m,t$ be some non-negative integers satisfying
$2k+2l+2m+t\leq n$. Set $\delta(0,0,0,0)=e$ and if at least one of
$k,l,m,t$ is not zero, set
\begin{multline}\label{eq:can}
\delta(k,l,m,t)=\theta_1\theta_3\cdots\theta_{2k-1}\xi_{2k+1}\xi_{2k+3}\cdots
\xi_{2k+2l-1}\nu_{2k+2l+1}\nu_{2k+2l+3}\cdots\cdot\\
\cdot\nu_{2k+2l+2m-1}\vartheta_{2k+2l+2m+1}
\vartheta_{2k+2l+2m+2}\cdots\vartheta_{2k+2l+2m+t}.
\end{multline} 
The element $\delta(k,l,m,t)$ such that $l=0$ or $m=0$ will be 
called a {\em canonical element} of  type $(k,l,m,n)$. 

\begin{corollary}\label{cor511}
Every $S_n\times S_n$-orbit contains a
canonical element.
\end{corollary}

\begin{proof}
Because of Corollary~\ref{cor530} we have to prove that,
the $S_n\times S_n$-orbit of the element
$\vartheta(A)\gamma_{i_1,j_1}\cdots\gamma_{i_s,j_s}$,
satisfying the conditions of Corollary \ref{cor530},
contains a canonical element.
Using conjugation, we can always reduce
$\vartheta(A)\gamma_{i_1,j_1}\cdots\gamma_{i_s,j_s}$
to some $\delta(k,l,m,t)$. However, it might happen that
both $m$ and $l$ are non-zero. Without loss of generality
we may assume $m\geq l\geq 1$. Using \eqref{pb4} and conjugation
we get that the $S_n\times S_n$-orbit of the element 
$\mu_{i,j}\nu_{p,q}$ contains $\epsilon_{i,j}\vartheta_p\vartheta_q$
provided that $\{i,j\}\cap\{p,q\}=\varnothing$. Hence the
$S_n\times S_n$-orbit of our $\delta(k,l,m,t)$ contains
$\delta(k+1,l-1,m-1,t+2)$. Proceeding by induction we get that
the $S_n\times S_n$-orbit of our $\delta(k,l,m,t)$ contains
$\delta(k+l,0,m-l,t+2l)$, which is canonical. This completes the proof.
\end{proof}

\begin{lemma}\label{lemma422}
The $S_n\times S_n$-orbits of the canonical element 
$\delta(k,l,0,t)$ and $\delta(k,0,l,t)$ contain at most 
\begin{displaymath}
\frac{(n!)^2}{(k+l)!2^{k+l}t!k!2^{k}(2l+t)!(n-2k-2l-t)!}
\end{displaymath}
elements.
\end{lemma}

\begin{proof}
We will prove the statement for the element $\delta(k,l,0,t)$. For
$\delta(k,0,l,t)$ the proof is analogous. We use the arguments similar 
to those from the proof of Lemma \ref{lemma222}. It is enough to show 
that the stabilizer of $\delta(k,l,0,t)$ under the $S_n\times S_n$-action
contains at least $(k+l)!2^{k+l}t!k!2^{k}(2l+t)!(n-2k-2l-t)!$ elements.  
Set 
\begin{gather*}\Sigma^0_i=\sigma_{2i}\sigma_{2i-1}\sigma_{2i+1}
\sigma_{2i},\,\, 1\leq i\leq k+l-1;\\
\Sigma^1_i=\sigma_{2i}\sigma_{2i-1}\sigma_{2i+1}\sigma_{2i}
\sigma_{2i-1},\,\,\,\, 1\leq i\leq k+l-1.
\end{gather*}
Then both $\Sigma^0_i$  and $\Sigma^1_i$ swap 
the sets $\{2i-1, 2i\}$ and $\{2i+1, 2i+2\}$. It follows that 
the group $H$, generated by all $\Sigma^0_i$, consists of all 
permutations of the set $\{1,2\}, \{3,4\},\dots, \{2k+2l-1, 2k+2l\}$ 
and is therefore isomorphic to the group $S_{k+l}$. It is further
easy to see that the group $\tilde{H}$, generated by all $\Sigma^0_i$ and 
$\Sigma^1_i$, is isomorphic to the wreath product $H\wr S_2$.
From \eqref{it51} and \eqref{it31} it follows that the left  
multiplications with $\Sigma^0_i$ and $\Sigma^1_i$ 
stabilizes $\delta(k,l,0,t)$. Therefore the left 
multiplication with each element of $\tilde{H}$ stabilizes 
$\delta(k,l,0,t)$ as well. 
Now, from \eqref{it43} and \eqref{pb2} it follows that 
\begin{equation*}
\sigma_i\eta_i=\sigma_i\vartheta_i\vartheta_{i+1}\theta_i=
\vartheta_{i+1}\sigma_i\vartheta_{i+1}\theta_i=
\vartheta_i\sigma_i\vartheta_i\theta_i=
\vartheta_i\vartheta_{i+1}\theta_i=\eta_i.
\end{equation*}
for all $i=1,\dots, n-1$. Moreover,
\begin{multline*}
\sigma_{i+1}\eta_i\eta_{i+2}=
\sigma_{i+1}\vartheta_{i+1}\theta_i\vartheta_{i+2}\theta_{i+2}=
\sigma_{i+1}\vartheta_{i+1}\vartheta_{i+2}\theta_i\theta_{i+2}=\\
\vartheta_{i+1}\vartheta_{i+2}\theta_i\theta_{i+2}=
\vartheta_{i+1}\theta_i\vartheta_{i+2}\theta_{i+2}=\eta_i\eta_{i+2}
\end{multline*}
for all $i=1,\dots, n-3$ by \eqref{pb1} and \eqref{it43} and
\begin{equation*}
\sigma_{i+1}\eta_i\vartheta_{i+2}=
\sigma_{i+1}\vartheta_{i+1}\theta_i\vartheta_{i+2}=
\sigma_{i+1}\vartheta_{i+1}
\vartheta_{i+2}\theta_i=\vartheta_{i+1}
\vartheta_{i+2}\theta_i=\eta_i\vartheta_{i+2}
\end{equation*}
for all $i=1,\dots, n-2$ again by \eqref{pb1} and \eqref{it43}.
Using this and the fact that  $\eta_i$ commutes with each of 
$\theta_j$, $\eta_j$, $\xi_j$ whenever $\vert i-j\vert> 1$ we see 
that each of the elements $\sigma_i$, $2k+2l-1\leq i\leq 2k+2l+t$,
stabilizes $\delta(k,l,0,t)$ under the left multiplication. 
All these elements generate the group $H_0\simeq S_{t}$, which 
stabilizes $\delta(k,l,0,t)$ and has trivial intersection with 
$\tilde{H}$. Let $H_1=H_0\times \tilde{H}$.

Analogously one shows that there is a group, $H_2$,
isomorphic to the wreath product 
$(S_{k}\wr S_2)\times S_{2l+t}$, such that 
each element of this group stabilizes $\delta(k,l,0,t)$ with 
respect to the right multiplication. Apart from this, from 
\eqref{it31} we have that conjugation by any element from 
the group $H_3=\langle \sigma_{2k+2l+t+1}, \dots, 
\sigma_{n-1}\rangle\simeq S_{n-2k-2l-t}$ stabilizes 
$\delta(k,l,0,t)$. Observe that the group, generated by 
$H_1$, $H_2$ and $H_3$, is a direct product of $H_1$, $H_2$ and $H_3$.
Hence, using the product rule we derive that 
the cardinality of the stabilizer of $\delta(k,l,0,t)$ is at least 
\begin{displaymath}
(k+l)!2^{k+l}t!k!2^{k}(2l+t)!(n-2k-2l-t)!,
\end{displaymath}
and the proof is complete.
\end{proof}

\begin{proof}[Proof of Theorem~\ref{theorem3}]
Comparing Lemma~\ref{lemma422} and
Proposition~\ref{comb}\eqref{comb.4} we have
$|T|\leq |\mathfrak{B}_n|$. Since $\varphi:T\to \mathfrak{B}_n$
is surjective we have $|T|\geq |\mathfrak{B}_n|$. Hence
$|T|=|\mathfrak{B}_n|$ and $\varphi$ is an isomorphism.
\end{proof}


\vspace{0.2cm}

\noindent
G.K.: Algebra, Department of Mathematics and Mechanics, Kyiv Taras
Shevchenko University, 64 Volodymyrska st., 01033 Kyiv, UKRAINE,
e-mail: {\tt akudr\symbol{64}univ.kiev.ua}
\vspace{0.2cm}

\noindent
V.M: Department of Mathematics, Uppsala University, Box. 480,
SE-75106, Uppsala, SWEDEN, email: {\tt mazor\symbol{64}math.uu.se}

\end{document}